\documentclass[a4paper]{amsart}

\usepackage[T1]{fontenc}
\usepackage{hyperref}
\usepackage{amsmath,mathtools}
\usepackage{amssymb}
\usepackage{amsthm}
\usepackage{graphicx}
\usepackage[bottom]{footmisc}
\usepackage[rgb]{xcolor}
\usepackage[margin=2cm]{geometry}
\usepackage{tikz}
\usepackage{subcaption}
\usepackage{float}
\usepackage{booktabs}

\usepackage{hyperref}
\hypersetup{
 colorlinks=true}
\usepackage{cite}

\usepackage{bm}

\DeclareMathAlphabet{\mathpzc}{OT1}{pzc}{m}{it}

\newcommand*\pFqskip{8mu}
\catcode`,\active
\newcommand*\pFq{\begingroup
\catcode`\,\active
\def ,{\mskip\pFqskip\relax}%
\dopFq
}
\catcode`\,12
\def\dopFq#1#2#3#4#5{%
    {}_{#1}F_{#2}\left(\genfrac..{0pt}{}{#3}{#4}\bigg\vert\,#5\right)%
\endgroup
}

\newtheorem{theorem}{Theorem}
\newtheorem{corollary}{Corollary}

\newtheorem{proposition}{Proposition}
\newtheorem{definition}{Definition}

\newcommand{\pf}{\noindent {\it Proof.}}
\renewcommand{\qed}{\hfill $\square$}

\DeclareMathOperator{\Ai}{Ai}

\newcommand\das{\ensuremath{\overset{d}{=}}}

\newcommand{\bigO}{\mathcal{O}}
\newcommand{\wbar}[1]{\overline{#1}}

\newcommand{\be}{\begin{equation}}
\newcommand{\ee}{\end{equation}}
\newcommand{\ef}[1]{\, #1}

\newcommand{\dd}{\mathrm{d}}

\newcommand{\lam}{\lambda}
\newcommand{\eps}{\epsilon}

\newcommand{\smfrac}[2]{\genfrac{}{}{0.25pt}{1}{#1}{#2}}

\newcommand{\mycirc}{\raisebox{2pt}{${\scriptstyle \circ}$}}

\begin{document}

\title{The $p$-Airy distribution}

\author[S. Caracciolo]{Sergio Caracciolo}
\author[V. Erba]{Vittorio Erba}  
\address[S. Caracciolo and V. Erba]{Dipartimento di Fisica, University
  of Milan and INFN, via Celoria 16, 20133 Milan, Italy}
\author[A. Sportiello]{Andrea Sportiello}
\address[A. Sportiello]{LIPN, and CNRS, Universit\'e Paris 13,
  Sorbonne Paris Cit\'e, 99 Av.~J.-B.~Cl\'ement, 93430 Villetaneuse,
  France}

\begin{abstract}
    In this manuscript we consider the set of Dyck paths equipped with
    the uniform measure, and we study the statistical properties of a
    deformation of the observable ``area below the Dyck path'' as the
    size $N$ of the path goes to infinity.  The deformation under
    analysis is apparently new: while usually the area is constructed
    as the sum of the heights of the steps of the Dyck path, here we
    regard it as the sum of the lengths of the connected horizontal
    slices under the path, and we deform it by applying to the lengths
    of the slices a positive regular function $\omega(\ell)$ such
    that $\omega(\ell) \sim \ell^p$ for large argument.  This shift
    of paradigm is motivated by applications to the Euclidean Random
    Assignment Problem in Random Combinatorial Optimization, and to
    Tree Hook Formulas in Algebraic Combinatorics.

    For $p \in \mathbb{R}^+ \smallsetminus \left\{ \frac{1}{2}\right\}$, we characterize
    the statistical properties of the deformed area as a function of
    the deformation function $\omega(\ell)$ by computing its integer
    moments, finding a generalization of a well-known recursion for
    the moments of the area-Airy distribution, due to Tak\'acs.  Most
    of the properties of the distribution of the deformed area are
    \emph{universal}, meaning that they depend on the deformation
    parameter $p$, but not on the microscopic details of the function
    $\omega(\ell)$.  We call \emph{$p$-Airy distribution} this family
    of universal distributions.

    Finally, we briefly study the limits $p\rightarrow\frac{1}{2}$,
    and $p\rightarrow 0$ (that is, with function $\omega(\ell) \sim
    \log(\ell)$), in which the recursion is singular in the sense of
    L'H\^{o}pital, and the determination of the moments is more
    subtle. A more detailed derivation of the singular cases will be
    developed in a companion paper.
\end{abstract}

\maketitle

\section{Introduction}
\label{sec.intro}

\noindent
The statistics of the area enclosed by a standard Brownian
excursion\footnote{I.e.\ a Brownian motion that starts at $(0,0)$,
  ends at $(1,0)$ and never reaches negative heights.} and the real
axis is well studied, and is governed by the so-called \emph{area-Airy
  distribution} $f_{\rm Ai}(x)$ (we follow the naming convention of
\cite{Flajolet}).  The area-Airy distribution owes its name to the
fact that its density and its Laplace transform admit a spectral
representation related to the zeros of the Airy function $\Ai(x)$
\cite{darling1983,takacs1991}.  An interesting fact, dating back to
the pionereeing work of Mark Kac \cite{Kac1951}, is that the
area-statistics of a Brownian excursion, a classical problem in
Probability Theory, is related to the study, in Quantum Mechanics, of
a one-dimensional particle subject to a linear potential and to a hard
wall in the origin, explaining the presence of the Airy function.

The area-Airy distribution has recently attracted attention in
Statistical Physics, where it appears to model a large number of
phenomena. A non-comprehensive list (more can be found in
\cite{Flajolet}) features the \emph{maximum height of fluctuating
  interfaces} \cite{Majumdar2005}, the \emph{size of avalanches in
  sandpile models} \cite{stapleton2006}, the \emph{size of ring
  polymers} \cite{medalion2016} and the \emph{anomalous diffusion in
  cold atoms} \cite{barkai2014}.  Very recently, the area-Airy
distribution function was measured experimentally for the first time
in a dilute colloidal system \cite{agranov2019}.  Generalizations of
the area-Airy distribution to the statistics of the area of other
Brownian processes, and to other properties of such Brownian
processes, can be found in
\cite{takacs1992,takacs1993,takacs1995,tolmatz2000,tolmatz2002,tolmatz2003,the2004,tolmatz2005}.

The present paper is aimed to the definition and study of a
one-parameter family of distributions, that we call 
\emph{$p$-Airy distributions}, which generalise the area-Airy case
(corresponding to $p=1$) to the range $p \in \mathbb{R}^+$.  In the
remaining of this introduction we give a list of probabilistic
problems in which our generalisation arises naturally.  As we will
see, similarly to the original Airy distribution, our
generalisation is ``universal'', that is, it does not depend on the
microscopic details of the system that leads to its definition, and,
for this reason, we expect that it can arise in a variety of
applications in Statistical Mechanics and Probability, on a similar
basis of the list of applications of the Airy distribution presented
above.  A \emph{different} generalization of the Airy distribution,
based on a phenomenon of coalescence of multiple saddle points, has
been pursued in
\cite{haug2017HigherOrderAiryScaling,haug2018MulticriticalScalingLattice}.
Our definition and list of examples passes through a \emph{d\'etour}
from stochastic processes to discrete combinatorics, along the line of
Donsker's theorem (in reverse), and not dissimilar in spirit from what
is done in \cite{takacs1991} for the case $p=1$.

\begin{figure}[t]
\[
\includegraphics[scale=.9]{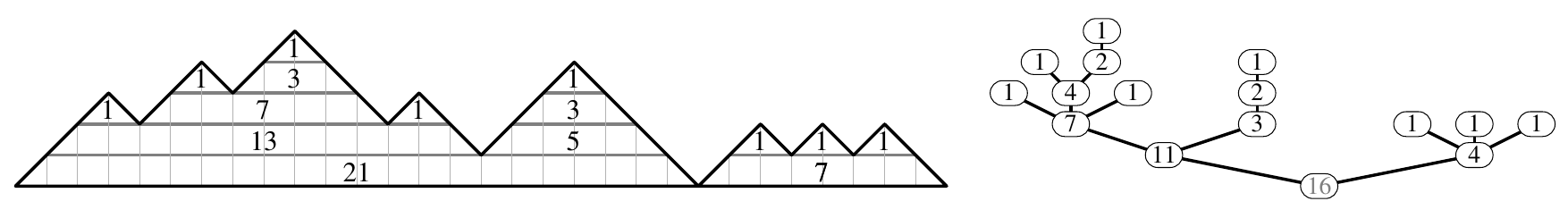}
\]
\caption{\label{fig.walktree}\footnotesize \textbf{Left:} a Dyck path of length $2N=30$,
  subdivided into connected horizontal slices. In the middle of each
  strip $e$, it is indicated the value of its length $\ell_w(e)$. 
  \textbf{Right:} the
  corresponding rooted planar tree $t(w)$. On each of the $N$ non-root
  vertices $v$, it is indicated the value of $h_t(v)$, while the root
  vertex $r$ has $h_t(r)=N+1$.}
\end{figure}

Brownian excursions are the continuum limit of a class of discrete
lattice paths called \emph{Dyck paths}.  A Dyck path of size $N$ is a
sequence of $N$ up- and $N$ down-steps, that is steps $(+1,+1)$ and
$(+1,-1)$ on the two-dimensional integer lattice, starting at
$(0,0)$, ending at $(2N,0)$ and never reaching negative heights.  The
area between a Dyck path $w$ and the real axis can be interpreted in
two natural ways, both as a sum of half-integers 
$\{h_w(i)\}_{1 \leq i \leq 2N}$ associated to the heights of the $2N$ vertical slices
of the walk (analogous to a Riemann-like integral approximation), or
as a sum of integers $\{\ell_w(e)\}_{1 \leq e \leq N}$ associated to the lengths of
the $N$ connected horizontal slices $\{e\}=E(w)$ of $w$ (analogous to
a Lebesgue-like integral approximation; see
Figure~\ref{fig.walktree}). That is,
\begin{equation}\label{eq:areaDyck}
A(w)=\sum_{i=1}^{2N} h_w(i) = \sum_{e \in w} \ell_w(e)
\ef,
\end{equation}
where with a slight abuse of notation we write $e \in w$ meaning 
$e \in E(w)$.  The area-Airy distribution must be recovered by taking
the continuum limit $N\rightarrow\infty$ of the distribution of $A(w)$
(induced by the uniform measure over Dyck paths), with a rescaling
factor $N^{\frac{3}{2}}$.  Indeed, also in \cite{takacs1991} the
characterization of the area-Airy distribution is based on the
continuum limit of the distribution of the area of Dyck paths.

When considering a stochastic quantity which is the sum of many local
contributions, it is natural to consider the distribution of the
\emph{moments} of the local terms. For example, the study of 
$A_p^{\rm BM}(w):=\sum_{i=1}^{2N} h_w(i)^p$, where the upper-script BM
refers to Brownian Motion, has a long tradition
\cite{mansuy2008aspects}. However, the analogous generalisation 
$A_p = \sum_{e \in w} \ell_w(e)^p$ appears to be new, and is indeed
the generalisation we are interested in.

One aspect of the aforementioned universality is that, in a sense that
we make precise later on, we can replace $\ell_w(e)^p$ by any function
$f(\ell)$ such that the large-$\ell$ behaviour is $f(\ell) \sim
\ell^p$. 
Thus, we define the generalization of Equation~\eqref{eq:areaDyck} as
\begin{equation}
\label{eq:areaDyckDeformed}
A_{(\omega_p)}(w) = \sum_{e \in w} 
\omega_p\!\left( \frac{\ell_{w}(e)-1}{2} \right)
\end{equation}
where $p\geq 0$, $w$ is a Dyck path, $( \ell_w(e) - 1)/2$ is the semi-length of the
horizontal slice $e$ in $w$ and $\omega_p$ is a positive regular
function on the integers such that 
\be
\label{eq.defeta}
\omega_p(k) \sim k^p \left( 1+\bigO\left( k^{-\eta} \right) \right)
\ee
for large $k$, and some $\eta > 0$.\;\footnote{The 
  subscript $p$ in the notation
  $\omega_p(k)$ is to stress the asymptotic behaviour of this
  function, and to prepare the notation for the case of interest in
  which $\{\omega_p(k)\}_{p \in \mathbb{R}^+}$ is a smooth family of functions.
  For the moment though, $p$
  is considered fixed and $\omega_p(k)$ is just one function of a
  single argument, $k$.} 

We want to compute the statistics of $A_{(\omega_p)}$ induced by the
uniform measure on two classes of lattice paths, namely the Dyck paths
that we have just described (or \emph{excursions} in the following)
and the Dyck bridges (which are Dyck paths without the non-negative
height constraint).\footnote{The names are slightly unusual in
  combinatorics, but are induced by the fact that the continuum limit
  of Dyck bridges and excursions, as we have called them, is the
  Brownian bridge and excursion, respectively.}  As a result, the
usual area of a Dyck path, or the unsigned area of a Dyck bridge, is
recovered with the choice
$\omega_p(k) = k+\frac{1}{2}$, up to the trivial overall factor of
$2$.

This generalization defines a deformation of the area-Airy
distribution and of its analogue for Brownian bridges, and it is
motivated by (at least) two rather different applications that we
summarize below.

\subsection{Euclidean Random Assignment Problem}

\begin{figure}[b!]
    \centering
    \begin{subfigure}{.49\textwidth}
\centering
\framebox{\includegraphics[width=4cm]{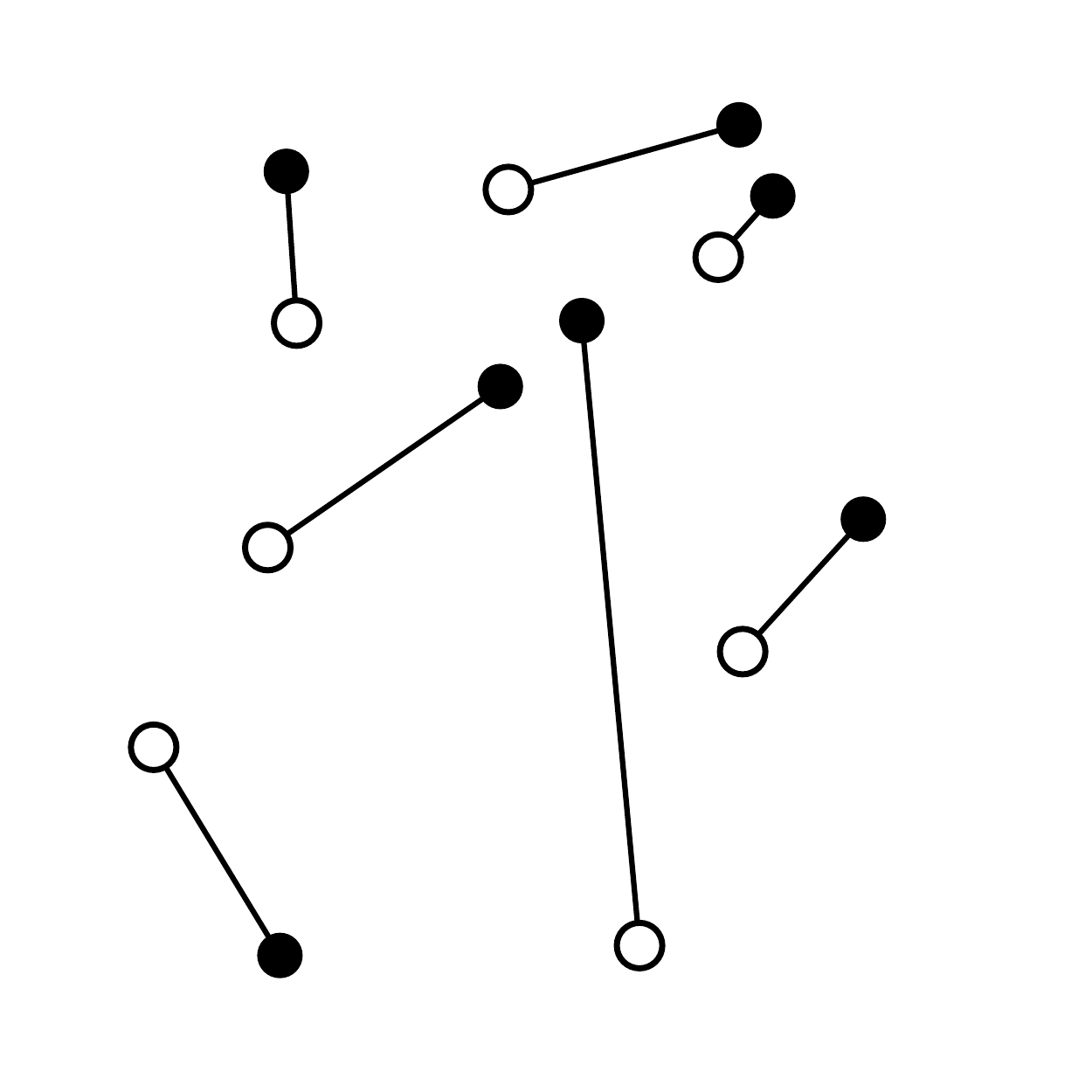}}
\label{fig:REAP1}
    \end{subfigure}%
    \begin{subfigure}{.49\textwidth}
\centering
\includegraphics[width=6cm,trim={0 3.5cm 0 3.5cm},clip]{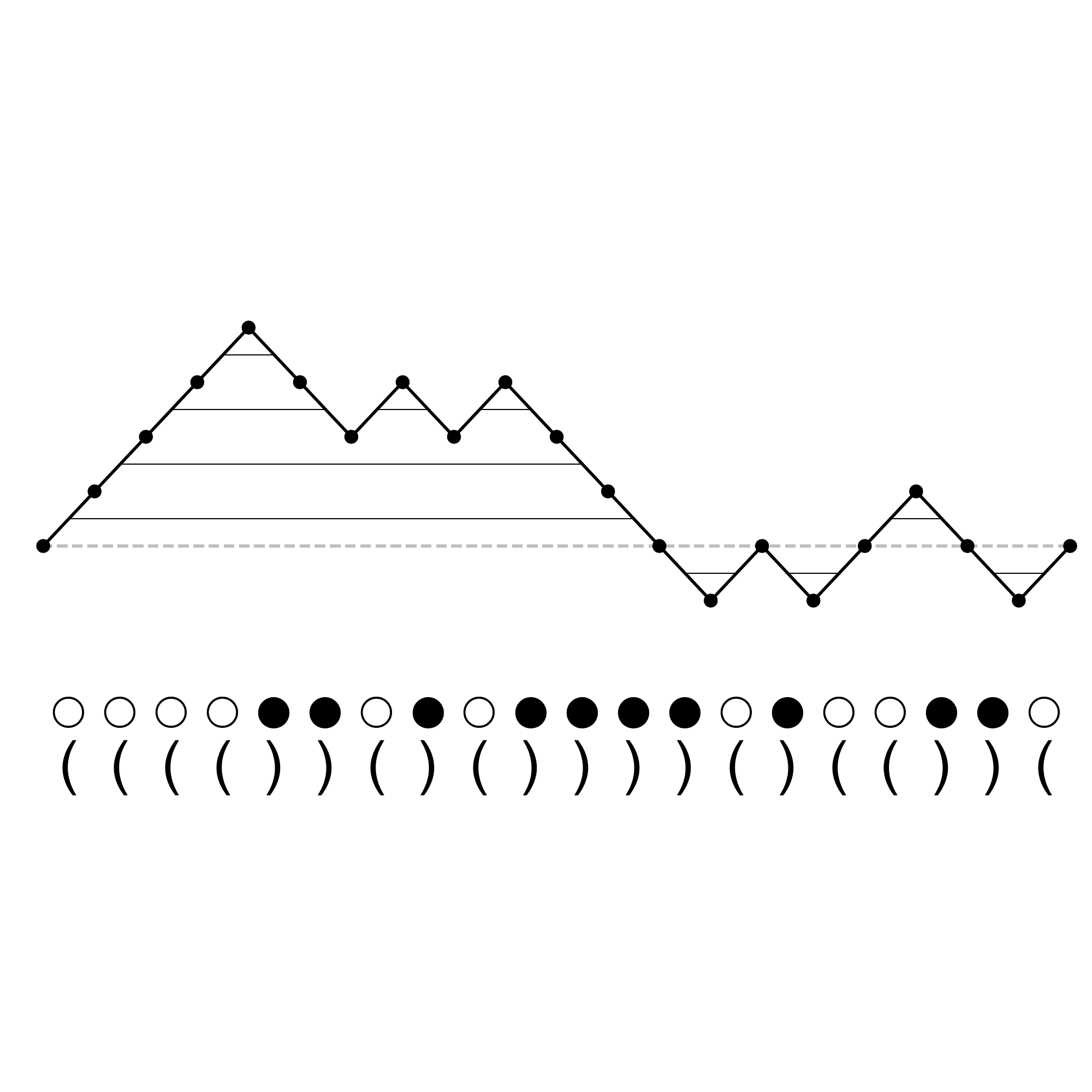}
\label{fig:REAP2}
    \end{subfigure}%
    \caption{\footnotesize 
\textbf{Left:} 
example of optimum matching with $N=7$ and $p=2$ in $[0,1]^{2}$.
\textbf{Right:} 
an example of the Dyck bridge of length $N=10$ associated to a sequence of white and black points.
Horizontal solid lines represent the Dyck matching associated to this configuration of points.
    }
    \label{fig:REAP}
\end{figure}

The Euclidean Assignment Problem (EAP) is a combinatorial optimization
problem in which one has to pair $N$ white points to $N$ black points
minimizing a certain cost function, depending on the Euclidean distance
among the points.
To be more precise, consider $N$ white points with coordinates
$\{w_{i}\}_{i=1}^{N}$
and $N$ black points with coordinates $\{b_{i}\}_{i=1}^{N}$,
and let $\tilde{c}(x)$ be a real function (in this work we will
consider the case $\tilde{c}(x)=|x|^{p}$). We call $\tilde{c}(x)$
the \emph{cost function} of the problem.

Solving one instance of the EAP corresponds to finding the perfect
matching\footnote{A \emph{perfect matching} is a bijection between the
  black and the white points, and a permutation $\pi$ describes the
  matching such that $w_{i}$ and $b_{j}$ are matched if and only if
  $\pi(i)=j$. See Figure~\ref{fig:REAP}.} between the white and
black points which minimizes the sum of the cost function of the
distances among the paired points, that is, a permutation 
$\pi \in \mathcal{S}_{N}$ such that the total cost
\begin{equation}
    \begin{split}
H[\pi] = \sum_{i=1}^{N} \tilde{c}(w_{i}-b_{\pi(i)}) \, 
    \end{split}
\end{equation}
is minimal.

The random version of the EAP (that we shall denote by ERAP) is the
probabilistic problem in which we study the measure over the instances
obtained by extracting the positions of the points according to some
probability law.  The total cost of the minimum matching becomes a
random variable and its statistical properties can be investigated
also as a function of the exponent $p$.  Plenty of results exist in
one dimension for the cases $p \geq
1$~\cite{Caracciolo:159,Caracciolo:160,Caracciolo:177}, also in a non
compact domain~\cite{Caracciolo:172} and some for $p<0$ in the
repulsive regime~\cite{Caracciolo:169}, due to the convexity of the
cost function $\tilde{c}$.  There have been remarkable extensions in
higher dimensions~\cite{Caracciolo:158,Caracciolo:162}, particularly
in the two-dimensional case when
$p=2$~\cite{Caracciolo:163,Ambrosio2016,Ambrosio2018,Ambrosio2019}.
Instead, for $0<p<1$, where the cost function is concave, the problem
is not equally-well understood.  It is known that the optimal matching
must satisfy certain nesting properties~\cite{mccann1999, DelonSS12}, but those
are not restrictive enough to fully characterize it.

Recently, two independent upper bounds to the average optimal cost for
$0<p<1$ were investigated by
combinatorial~\cite{Caracciolo:159,Caracciolo:180} and
measure-theoretic \cite{Bobkov2019} means.  In particular,
in~\cite{Caracciolo:180} the authors studied the cost of the
\emph{Dyck matching} (in the following, the \emph{Dyck cost}), and
computed the asymptotic properties of its average value in the limit
of $N\rightarrow\infty$.  The Dyck matching associated to a certain
configuration of points is determined by the $N$ horizontal slices
$\{e\}=E(w)$ of the Dyck bridge $w$ generated by scanning over the
points ordered by increasing coordinate, and performing an up
(resp.\ down) step for each white (resp.\ black) point encountered
(see Figure~\ref{fig:REAP}).  This relation between matchings and
lattice paths has already been exploited in the statistical mechanics
literature to study, for example, the secondary structure of folded
RNA in the models considered 
in~\cite{Tamm_2007,hofacker1998combinatorics}.

In the case in which the points are equispaced, the distance between
two matched points in the Dyck matching and the length of the
corresponding horizontal slice of the Dyck bridge are equivalent, so
that the Dyck cost is exactly given by
\begin{equation}
    \begin{split}
H^{(p)}_{\rm Dyck}(w) 
&= \sum_{e \in w} (\ell_{w}(e))^{p} = 2^{p}
A_{\left(\omega_p\right)}(w) \quad \text{with } \quad \omega_p(k) =
\left( k+\frac{1}{2} \right)^p 
\, ,
    \end{split}
\end{equation}
where $w$ is the Dyck bridge associated to the position of the points.
Thus, studying the statistics of the Dyck cost in a model of
equispaced points is a special case of our original problem of
determining the distribution of $A_{(\omega_p)}(w)$ induced by the
uniform measure over Dyck bridges of length~$N$.

We conclude this section by pointing out that in \cite{Caracciolo:180}
the authors highlighted a powerful universality for the statistics of
the Dyck cost.  If $p\geq\frac{1}{2}$, the average Dyck cost is vastly
independent on the model of spacings between the points; in other
words, it depends only on the long-distance behaviour of the cost
function $\tilde{c}(k) \sim k^p$, and not on its microscopic details
at short distance. We will find again this property, now at the level
of the distribution of this stochastic quantity, not only of the
average.

\subsection{Tree Hook Formula, and the statistics of star subtrees}
\label{ssec.hookform}

It is well known that Dyck paths of length $2N$ and rooted planar
trees with $N+1$ vertices are in bijection, and this bijection relates
the non-root vertices $v$ of the tree and the horizontal slices $e$ of
the path (see e.g.\ \cite[sec.\ I.5]{flajolet2009}, and Figure
\ref{fig.walktree}).

Let $w$ be a path of length $2N$, and $t=t(w)$ the corresponding tree.
The \emph{hook} $h_t(v)$ at a vertex $v \in V(t)$ is defined as the
number of vertices in the sub-tree of $t$ rooted at $v$ (including
$v$). Under the bijection we just have that, if the (non-root) vertex
$v$ of the tree is in bijection with the edge $e$, then
\begin{equation}
\label{eq.lh1}
\ell_w(e)=2h_t(v)-1
\ef,
\end{equation}
while obviously the root vertex $r$ has $h_t(r)=N+1$, for all trees
of size $N+1$.
The product of the hooks enters the celebrated ``tree hook formula''
\cite[sec.\ 5.1.4, ex.\ 20]{knuthACP3} (see e.g.\ \cite{FerGou12} for
an introduction to the tree-hook formula and its generalisations),
which counts the fraction $L(t)$ of labellings of the vertices of a
rooted tree which are \emph{increasing}:
\begin{equation}
\label{eq.thf}
L(t) =
\prod_v h_t(v)^{-1} = \exp
\bigg[ -\sum_v \ln h_t(v) \bigg] = 
\lim_{p \to 0}
\exp\bigg[-\frac{1}{p} 
\bigg(
\sum_v h_t(v)^p - (N+1)
\bigg)\bigg]
\ef.
\end{equation}
So, the distribution of $L(t)$ induced by the uniform distribution on
rooted planar trees with $N+1$ vertices is related to the statistics
of the quantity 
\begin{equation}
    \begin{split}
H^{(p)}_{\rm tree}(t) = \sum_v (h_t(v))^p = A_{\left(\omega_p
  \right)}(w(t)) + (N+1)^{p}  
\qquad \text{with} \quad \omega_p(k) = (k+1)^p      \, 
    \end{split}
\end{equation}
in the limit $p \to 0$.

Moreover, the quantity $H_{\rm tree}^{(p)}(t)$, for $p$ integer,
counts the $(p+1)$-tuples $(v_0,v_1,\ldots,v_p)$ of vertices of $t$
such that $v_0 \preceq v_i$ (i.e.\ $v_0$ is on the path connecting
$v_{i}$ to the root) for all $i=1,\ldots,p$, while the analogous
quantity with $(h_t(v))^p$ replaced by $(h_t(v)-1)^p$ counts the
$(p+1)$-tuples as above, with $\preceq$ replaced by $\prec$.

\begin{figure}[tb!]
\[
\includegraphics[scale=.9]{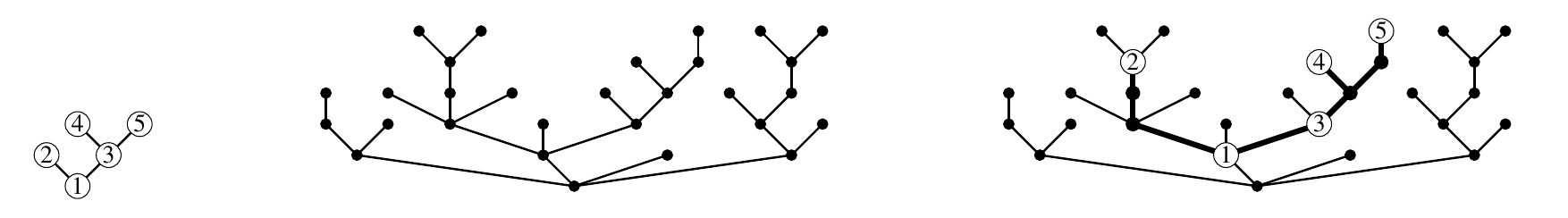}
\]
\caption{\footnotesize 
\textbf{Left:} a planar rooted labeled tree $u$.
\textbf{Center:} a planar rooted tree $t$.
\textbf{Right:} a proper embedding of $u$ into $t$.
}
\label{fig:treemb}
\end{figure}

This is a special case of the statistics $H(t,u)$, which, for $t$ and
$u$ two rooted trees, counts the number of \emph{embeddings} of $u$ in
$t$ (or \emph{proper embeddings}, for the $\prec$ case), where a
(proper) embedding of a rooted tree $u$ in the rooted tree $t$ is a
map from the vertices of $u$ into those of $t$ that preserves the
order relation $\preceq$ (or $\prec$), see an example in Figure
\ref{fig:treemb}.  The statistics of the number of embeddings of a
given rooted tree $u$, of size $\bigO(1)$, into random uniform trees
$t$ taken uniformly among all planar rooted trees of size $N$, in the
limit of large $N$, is a problem of separate interest. In the case in
which $u$ is the `star tree', composed of a root connected to $p$
children, this is related to the statistics of $A_{(\omega_p)}(w)$.
In this case, $w$ should be uniformly distributed over Dyck paths of
length $N$, and $\omega_p(k)=(k+1)^p$ or $k^p$ in the case of
embeddings or proper embeddings, respectively.

\subsection{A glance to our results}
\label{ssec.results}

In this manuscript, we fully characterize the distribution of the
random variable $A^{\rm (E/B)}_{(\omega_p)}(w)$ induced by the uniform
distribution over Dyck excursions or bridges for a generic cost
function $\omega_p$ as above, in the limit of $N\rightarrow\infty$.

The statistical ensembles of Dyck excursions and Dyck bridges are
treated with minor differences, and the results are analogous.  So,
for sake of compactness, we state here the main results for the case
of excursions (and drop the `E' superscript), and refer the reader to
the following Sections for the details in the case of Dyck bridges.
As anticipated by some facts already presented in
\cite{Caracciolo:180}, we identify three different behaviours as the
parameter $p$ varies.

For $p>\frac{1}{2}$, we find that
\begin{equation}
\label{eq:introPmaggiore}
    \begin{split}
A_{(\omega_p)} \das x_p N^{p+\frac{1}{2}} 
\big( 1 + 
\bigO\big( N^{-\min\left( p-\frac{1}{2}, \eta,\frac{1}{2} \right) } \big)
\big)
    \end{split}
\end{equation}
where `$\das$' means `distributed as', $\eta$ controls the error term
in the asymptotic behaviour of $\omega_p(k)$ (as in
Equation~(\ref{eq.defeta})) and $x_p$ is a random variable whose
integer moments are given by
\begin{equation}
\label{eq:introFattorelli}
    \begin{split}
\left\langle 
x_p^s
\right\rangle 
= \frac{4 \sqrt{\pi} s!}{\Gamma\left( s\left( p+\frac{1}{2}
  \right) -\frac{1}{2} \right)} 
\mu_s(p)
\, 
    \end{split}
\end{equation}
where
the coefficients 
$\mu_s(p)$ 
satisfy the quadratic recursion
\begin{equation}
\label{eq:recMuEintro}
    \begin{split}
\mu_{0}(p) 
&= -\frac{1}{2}
\\
\mu_{s}(p) 
&= 
\mu_{s-1}(p) 
\frac{\Gamma\left(
  s\left( p+\frac{1}{2} \right) -1  \right)}{2\Gamma\left(
  s\left( p+\frac{1}{2} \right) -1-p\right)}     +
\sum_{k=1}^{s-1} 
\mu_{k}(p)\mu_{s-k}(p) 
\qquad
\text{for} \quad s\geq1
\, .
    \end{split}
\end{equation}
This sequence of moments defines a unique distribution on
$[0,+\infty)$, 
that we call $\rho^{\rm (E)}(x;p)$, (where `(E)' stands for
`excursions'), and for $p=1$ reduces to a known recursion for the
moments of the area-Airy distribution $\rho^{\rm (E)}(x;1)$, due to
Tak\'acs.  We recall some interesting facts on the area-Airy
distribution in Appendix~\ref{app:airyfacts}. In particular, the
Tak\'acs recursion for the moments is given in
Equation~\eqref{eq:recTakacs}.  We stress again that the behaviour of
$A_{(\omega_p)}$ in this regime is \emph{universal}, in the sense that
it depends only on the asymptotic behaviour of the cost function, but
not on its details.  For bridges, an analogous result holds, with a
different moment recursion, see Proposition~\ref{prop:EBleading} later
on. This leads to the definition of a different family of
distributions
$\rho^{\rm (B)}(x;p)$.

For $0<p<\frac{1}{2}$, we find that the distribution of
$A_{(\omega_p)}(w)$ peaks around its average value 
$\alpha(\omega_p) N$, where $\alpha(\omega_p)$ is a non-universal
constant that depends on the details of $\omega_p$ that will be
defined in detail in Equation~\eqref{eq:LhatE0}, and which we
anticipate here to be determined by
\begin{equation}
\label{eq:LhatE0antic}
    \begin{split}
\sum_{k \geq 0}
\frac{\Gamma(k+\frac{1}{2})}{2\sqrt{\pi}\,\Gamma(k+2)}
\omega_p(k) z^k
=
\alpha(\omega_p)
+ \frac{\Gamma(p-\frac{1}{2})}{2 \sqrt{\pi}}
(1-z)^{\frac{1}{2}-p}
+ \cdots
    \end{split}
\end{equation}
where the dots stand for higher powers of $(1-z)$.

Nonetheless, the typical fluctuations around the mean are again of
order $N^{p+\frac{1}{2}}$, and universal, and their distribution is
given by the family of $\rho^{\rm (E)}(x;p)$, which is determined by
the same formulas (\ref{eq:introFattorelli}) and
(\ref{eq:recMuEintro}), here for $0<p<\frac{1}{2}$. An important
difference is that the distributions $\rho^{\rm (E)}(x;p)$ for
$p<\frac{1}{2}$ have support on $\mathbb{R}$, contrarily to the case
$p>\frac{1}{2}$, where the support is $\mathbb{R}^+$.  For bridges we
have the same picture, with average value $\alpha(\omega_p) N$ (the
same as for the case of excursions), and with the distribution of
fluctuations given by $\rho^{\rm (B)}(x;p)$.

These results can be summarized by saying that, for $p\neq\frac{1}{2}$
\begin{equation}
\label{eq:decomposition1}
\begin{split}
A^{\rm (E/B)}_{(\omega_p)} \das \alpha(\omega_p) N + 
x^{\rm (E/B)}_p N^{p+\frac{1}{2}}
\big( 1 + 
\bigO\big( N^{-\min\left( p, \eta,\frac{1}{2} \right) } \big)
\big)
\, ,
\end{split}
\end{equation}
with $x^{\rm (E/B)}_p$ distributed with law $\rho^{\rm (E)}(x;p)$ for
excursions, and with law $\rho^{\rm (B)}(x;p)$ for bridges.  Note
that, even just for $p>\frac{1}{2}$, this claim is slightly stronger
than (\ref{eq:introPmaggiore}), as we have explicitated the possible
correction term scaling as $N^{p-\frac{1}{2}}$.  The proof of these
claims is detailed in Section~\ref{sec:integerMoments}, for what
concerns the explicit formulas, and in Section~\ref{sec:uniqueness},
for what concerns the uniqueness.

The integer moments of $\rho^{\rm (E)}(x;p)$ have a nice combinatorial
interpretation in terms of a sum over rooted planar trees, weighted in
a peculiar way. We prove this connection in
Section~\ref{sec:combinatorial}. For what we know, this fact, which
remains non-trivial also for the $p=1$ ordinary Airy distribution, and
which, in this case, could have been evinced also from Tak\'acs
recursion, was not previously observed.

If $\{\omega_p(k)\}_{p \in \mathbb{R}^+}$ is a family of functions that depend smoothly on
$p$, we can perform the limit $p \to \frac{1}{2}$, but the
involved procedure is delicate.  In fact, all the moments of
$\rho^{\rm (E)}(x;p)$ diverge as $p$ tends to $\frac{1}{2}$, and so
does $\alpha(\omega_p)$.
We anticipate the fact, to be proven in a forthcoming companion paper,
that there exists a family of functions $t(p)$ such that all the
moments of $\rho^{\rm (E)}(x;p)$ are simultaneously regularized by
shifting our random variable, $x \rightarrow x - t(p)$, and also such
that $\alpha(\omega_p) + t(p)$ is a regular function of $p$.  The
family of $t(p)$'s is an affine space and has the following
characterisation: $t$ can be extended to a meromorphic function in
$\mathbb{C}$, and it has a unique simple pole on the real-positive
axis at $p=\frac{1}{2}$, with coefficient
\be
\label{eq.deftstar}
t^* := \frac{1}{2\sqrt{\pi}}
\ee
(see Equation~\eqref{eq:tstar}).  We discuss these claims in
Section~\ref{sec:p05}.

In light of this claim on the $p=\frac{1}{2}$ case, when
$\{\omega_p(k)\}_{p \in \mathbb{R}^+}$ is a family of functions that depend smoothly on $p$,
Equation~\eqref{eq:decomposition1} is better rewritten as
\begin{equation}
\label{eq:decomposition2}
    \begin{split}
A_{(\omega_p)}^{\rm (E/B)}
 &\das \left( \alpha(\omega_p) + t(p)
N^{p-\frac{1}{2}} \right) N 
+ \tilde{x}^{\rm (E/B)}_p N^{p+\frac{1}{2}} 
\big( 1 + 
\bigO\big( N^{-\min\left( p, \eta,\frac{1}{2} \right) } \big)
\big)
\, ,
    \end{split}
\end{equation}
where $\tilde{x}^{\rm (E/B)}_p$ is stochastic, distributed with the
shifted law $\tilde{\rho}^{\rm (E/B)}(x;p)$ (implicitly, depending on
$t(p)$) such that
\begin{equation}
    \begin{split}
\tilde{\rho}^{\rm (E/B)}(x;p) 
= \rho^{\rm (E/B)}(x-t(p);p) \,,
    \end{split}
\end{equation}
and $\alpha(\omega_p)$ and $t(p)$ are deterministic constants (for a
given $p$ and $\omega_p$). Analogously, we will call
$\tilde{\mu}_s(p)$ the (purportedly finite) shifted moments.

The limit $p\rightarrow\frac{1}{2}$ of
Equation~\eqref{eq:decomposition2} is slightly non-trivial, as it
involves the use of L'H\^opital's rule, and gives
\begin{equation}
\label{eq:decomposition3}
    \begin{split}
A_{(\omega_{\frac{1}{2}})} 
&\das t^{*}\, N \log N 
+ \left( \alpha_0 + t_0 + \tilde{x}_{\frac{1}{2}} \right) N  + o(N)\,,
    \end{split}
\end{equation}
where $t^* = 1/(2\sqrt{\pi})$ (as defined in (\ref{eq.deftstar})),
and 
$\alpha_0$,
$t_0$ are the constant terms of the
Laurent
expansion at $p=\frac{1}{2}$ of the corresponding quantities,
i.e.\ 
\begin{align}
t(p) &= \frac{t^*}{p-\frac{1}{2}} + t_0 + \bigO(p-\smfrac{1}{2})
&
\alpha(\omega_p) 
&= -\frac{t^*}{p-\frac{1}{2}}
+\alpha_0 + o(1)
\,.
\end{align}
Using this representation, it is easy to see that at $p=\frac{1}{2}$
the distribution $A_{(\omega_p)}$ concentrates around its average
value $t^* N \log N$ with typical fluctuations of order $N$.  The
fluctuations are distributed with universal law 
$\tilde{\rho}^{\rm (E/B)}\left(x;\frac{1}{2}\right)$, apart from a
non-universal shift $\alpha_0 + t_0$ that depends on the details of
the cost function $\omega_p$ and on the regularization.  
Expressions for
the shifted moments $\tilde{\mu}_s(\frac{1}{2})$ are given, without
proof, in Section \ref{sec:p05}.

For the sake of concreteness, it is tempting to choose once and for
all a `canonical' regularization $t(p)$, within the family of possible
choices.  A natural choice is to shift the distributions in order to
have zero mean for each value of $p$, which corresponds to the choice
\begin{equation}
\label{eq.tpCano}
    \begin{split}
t(p) = \frac{ \Gamma(p-\frac{1}{2})}{2\,\Gamma(p)}
=
\frac{1}{2 \sqrt{\pi} \left(p-\frac{1}{2}\right)}
+\frac{\ln 2}{\sqrt{\pi}}
+\frac{p-\smfrac{1}{2}}{\sqrt{\pi}}
  \left((\ln 2)^2
-\frac{\pi^{2}}{12}\right)
+\bigO \big(\big(p-\smfrac{1}{2}\big)^2\big)
\, .
    \end{split}
\end{equation}
In this paper we will call the distributions 
$\tilde{\rho}^{\rm (E/B)}(x;p)$ resulting from the shift
(\ref{eq.tpCano}) `the'
\emph{$p$-Airy distribution} (for excursions and bridges
respectively), leaving aside the forementioned arbitrariness of the
regularisation.  We stress again that, under this choice,
$\tilde{\rho}^{\rm (E)}(x;1)$ is not the ordinary area-Airy
distribution, but rather its centered version, that is, the
distribution shifted by $-\sqrt{\pi}/2$.

The limit $p \to 0$ is trivial, especially under the choice
$\omega_0(k)=1$, that gives at sight in (\ref{eq:areaDyckDeformed})
$A_{(\omega_0)}(w) = N$ for all $w$. We just get $\alpha(\omega_0)=1$
and $\rho^{\rm (E)}(x;0)=\delta(x)$. Nonetheless, given a family of
functions $\{\omega_p\}_{p \in \mathbb{R}^+}$, one can study a more
subtle limit of the quantity $A_{(\omega_p)}(w)$, namely
\begin{equation}
    \begin{split}
\lim_{p\rightarrow 0^+} \frac{A_{(\omega_p)}(w) - N}{p} \, 
    \end{split}
\end{equation}
which amounts to studying the observable $A_{(\omega_{0^+})}$, where
$\omega_{0^+}$ denotes a function such that $\omega_{0^+}(k) \sim \log
k$ for large argument $k$.  We prove that, in this case, the
distributions peak around their average value in the limit of 
$N \rightarrow \infty$, with Gaussian fluctuations:
\begin{equation}
    \begin{split}
A_{(\omega_{0^+})} 
\das \alpha(\omega_{0^+}) N 
+ \sqrt{\gamma_E N \log N} \, z 
+ \bigO\big( \sqrt{N} \big)  \,
    \end{split}
\end{equation}
where $\gamma_E$ is the Euler's gamma constant, and $z$ is a standard
Gaussian random variable. We present these results in Section~\ref{sec:p0log}.


\section{The distribution of $A_{(\omega_p)}$ for $p\neq\frac{1}{2}$}
\label{sec:proofNo05}

\noindent
Our strategy to prove
Equation~\eqref{eq:decomposition1}, namely
\begin{equation}\label{eq:decomposition1sec1}
    \begin{split}
        A_{(\omega_p)} \das \alpha(\omega_p) N + x_p N^{p+\frac{1}{2}}\left(1 +
            \bigO\big(N^{-\min(p,\eta,\frac{1}{2})}\big) \right) \, ,
    \end{split}
\end{equation}
goes through the derivation of the asymptotic behaviour of the integer moments of $A_{(\omega_p)}$.
This will be achieved in Section~\ref{sec:integerMoments} by studying their generating functions, and then
applying singularity analysis.
The identification of the distribution, by the classical Carleman's
condition, will follow in Section~\ref{sec:uniqueness}.

\subsection{The integer moments of  $A_{(\omega_p,\epsilon)}$}
\label{sec:integerMoments}

As a first step, it is useful to consider a slight generalization of
$A_{(\omega_p)}(w)$, that is:
\begin{equation}
\label{eq:areaDyckDeformedShift}
A_{(\omega_p, \eps)}(w) 
= \sum_{e \in w} \left[ \omega_p \left( \frac{\ell_{w}(e)-1}{2} \right) - \eps \right]
= A_{(\omega_p)} - \eps N \, .
\end{equation}
The idea is that, in the following, $\eps$ can be tuned to cancel
exactly the deterministic term $\alpha(\omega_p)N$ in
Equation~\eqref{eq:decomposition1sec1}, that dominates when
$0<p<\frac{1}{2}$, thus allowing for a unified treatment of the two
regimes $0<p<\frac{1}{2}$ and $p>\frac{1}{2}$.

We shall consider the integer moments of $A_{(\omega_p,\eps)}(w)$ defined as
\begin{align}
    M^{\rm (E)}_{s}(N;\omega_p,\eps) &= C_{N}^{-1}\sum_{w \in \mathcal{C}_{N}}
 (A_{(\omega_p,\eps)}(w))^{s} \, , \\
    M^{\rm (B)}_{s}(N;\omega_p,\eps) &= B_{N}^{-1}\sum_{w \in \mathcal{B}_{N}}
 (A_{(\omega_p,\eps)}(w))^{s} \, ,
\end{align}
where $\mathcal{C}_{N}$ and $\mathcal{B}_{N}$ are the sets of Dyck
paths and Dyck bridges, and $C_{N} = |\mathcal{C}_N| =
\frac{1}{N+1}\binom{2N}{N}$ and $B_{N}=|\mathcal{B}_N|=\binom{2N}{N}$
denote the respective cardinalities.

To compute these moments, we will proceed as follows:
\begin{enumerate}
    \item we introduce two generating functions, $E_s(z)$ and
      $B_s(z)$, for the moments $M_s^{\rm (E/B)}(N;\omega_p,\eps)$,
      see Equation~\eqref{eq:Ez} and Equation~\eqref{eq:Bz}.  The
      analysis of the leading singular behaviour of $E_s(z)$ and
      $B_s(z)$ around their dominant pole determines the asymptotic
      behaviour for large $N$ of $M^{\rm(E/B)}(N;\omega_p,\eps)$.
    \item we show in Proposition~\ref{prop:recursionGF} how the
      generating functions $E_s(z)$ and $B_s(z)$ can be computed
      recursively, in terms of a suitable linear operator $\hat{L}$
      acting on the space of formal power series, and determined
      by~$\omega_p$;
    \item in the perspective of the asymptotic analysis mentioned
      before, we provide some tools, namely Proposition~\ref{prop:Lq}
      and Proposition~\ref{prop:epsilon}, to study the recursion
      relations of Proposition~\ref{prop:recursionGF} at their leading
      singular order.
\end{enumerate}
Summing all these efforts up, we will be able to prove
Proposition~\ref{prop:EBleading}, where we determine the coefficient
of the leading singular behaviour of $E_s(z)$ and $B_s(z)$, obtaining
that the integer moments of the random variable $x_p$ satisfy
the Equations~\eqref{eq:introFattorelli} and
\eqref{eq:recMuEintro} already presented in the Introduction
(and the respective versions for bridges).

\vspace{0.5cm}

\noindent
So, we start by introducing the generating functions (the dependence on $\omega_p$ and $\eps$ is understood)
\begin{align}
    \begin{split}\label{eq:Ez}
E_{s}(z) &= \frac{1}{2} \sum_{N\geq0} \frac{C_{N}}{s!} M_s^{\rm (E)}(N;\omega_p,\eps) \left( \frac{z}{4} \right)^{N}
= \frac{1}{2} \sum_{N\geq0} \sum_{w \in \mathcal{C}_{N}}\frac{(A_{(\omega_p,\eps)}(w))^{s}}{s!} \left( \frac{z}{4} \right)^{N} \, ,
    \end{split}
\\
    \begin{split}\label{eq:Bz}
B_{s}(z) &= \sum_{N\geq0} \frac{B_{N}}{s!} M_s^{\rm (B)}(N;\omega_p,\eps) \left( \frac{z}{4} \right)^{N}
= \sum_{N\geq0} \sum_{w \in \mathcal{B}_{N}}\frac{(A_{(\omega_p,\eps)}(w)^{s}}{s!} \left( \frac{z}{4} \right)^{N} \, .
    \end{split}
\end{align}
Let $\odot$ denote the Hadamard product between formal power series
(see e.g.\ \cite[V.3.2]{flajolet2009})
\begin{equation}
    \begin{split}
\bigg( \sum_{N\geq0} g_{N} z^{N} \bigg)  \odot 
\bigg( \sum_{N\geq0} h_{N} z^{N} \bigg) 
= \sum_{N\geq0} g_{N}h_{N}z^{N} \, 
    \end{split}
\end{equation}
and let
\begin{equation}
    L(z;\omega_p,\eps) = \sum_{k \geq 0} 
\left[\, \omega_p \left( k \right) - \eps \,\right] z^k \, .
\end{equation}
We will need the following definition:
\begin{definition}
Let $\hat{L}_{(\omega_p,\eps)}$ be the linear operator on formal power
series defined by
\begin{equation}
 \begin{split}
     \hat{L}_{(\omega_p,\eps)}[f](z) 
     &= L(z;\omega_p,\eps) \odot f(z) 
     = \sum_{N\geq0} f_{N} 
\left[\, \omega_p \left( N \right) - \eps \,\right] z^{N} 
     = \sum_{N\geq0} f_{N} \omega_p \left( N \right) z^{N} - \eps f(z) \, .
\end{split}
\end{equation}
\end{definition}
\noindent
Two small remarks are in order.
First, we can define the integer powers $\hat{L}_{(\omega_p,\eps)}^k$ as
\begin{equation}
 \begin{split}
     \hat{L}_{(\omega_p,\eps)}^k[f](z) 
     &= 
L(z;\omega_p,\eps)^{\odot k} \odot f(z)
= \sum_{N\geq0} f_{N} 
\left[\, \omega_p \left( N \right) - \eps \,\right]^k z^{N} \, .
\end{split}
\end{equation}
Furthermore, calling
$\hat{L}_{\omega} = \hat{L}_{(\omega,0)}$, we have
$\hat{L}_{(\omega_p,\eps)} = \hat{L}_{\omega_p} - \eps \mathbb{I}$, 
where $\mathbb{I}$ is the identity operator on generating
functions, that can be represented either by usual multiplication by
1, or by Hadamard multiplication by $(1-z)^{-1}$.

\begin{proposition}\label{prop:recursionGF}
    The generating functions $E_{s}(z)$ and  $B_{s}(z)$ satisfy, respectively,
    \begin{equation}\label{eq:recE}
\begin{split}
E_{0}(z) &= \frac{1-\sqrt{1-z}}{z} \, ,\\
E_{s}(z) &= \frac{z}{2\sqrt{1-z}}
\sum_{\substack{s_{1},s_{2},s_{3}\geq0
\\ s_{1}+s_{2}+s_{3}=s\\s_{1},s_{2}<s}} E_{s_{1}}(z)
\frac{1}{s_{3}!} \hat{L}_{(\omega_p,\eps)}^{s_3}[E_{s_{2}}](z) \qquad {\rm
for }\quad s\geq 1 
\,.
\end{split}
    \end{equation}
    and
    \begin{equation}\label{eq:recB}
\begin{split}
B_{0}(z) &= \frac{1}{\sqrt{1-z}} \, ,\\
B_{s}(z) &= \frac{z}{\sqrt{1-z}}  \sum_{\substack{s_{1},s_{2},s_{3}\geq0 \\ s_{1}+s_{2}+s_{3}=s\\s_{1}<s}} B_{s_{1}}(z) \frac{1}{s_{3}!} \hat{L}^{s_3}_{(\omega_p,\eps)}[E_{s_{2}}](z)
\qquad {\rm for }\quad s\geq 1 \, . 
\end{split}
    \end{equation}
    Notice that the range of the sums of Equation~\eqref{eq:recE} and
    Equation~\eqref{eq:recB} are slightly different: in the second case,
    $s_{2}=s$ is allowed.
\end{proposition}
\noindent
The proof can be found in Appendix~\ref{proof:recursionGF}.  For
excursions, the idea is to decompose a Dyck path $w$ as
$w=(+,w_{1},-,w_{2})$, where, for some $0\leq m \leq N-1$,
$w_{1} \in \mathcal{C}_{m}$, $w_{2} \in \mathcal{C}_{N-m-1}$ and $+/-$
are up/down step.  This allows to decompose
$A_{(\omega_p,\epsilon)}(w)$ into two independent terms relative to
$w_{1}$ and $w_{2}$ respectively.  A similar decomposition holds for
Dyck bridges, leading to the second set of Equations. Note that
(\ref{eq:recE}) is quadratic in the $E_s$'s, while (\ref{eq:recB}) is
linear in the $E_s$'s and $B_s$'s, so that a way to proceed is first
solve the non-linear recursive relation (\ref{eq:recE}), and then
deduce an inhomogeneous linear equation for the $B_s$'s from
(\ref{eq:recB}). This is one of the reasons why we mostly concentrate
on the case of excursions, and only sketch the minor modifications
required for bridges, when there are no substantial differences.

The exact treatment of Equation~\eqref{eq:recE} and
Equation~\eqref{eq:recB} and the expansion of $E_{s}(z)$ and
$B_{s}(z)$ are not possible in general, mainly due to the fact that
the operator $\hat{L}$ is complicated.  Nonetheless, we are able
to obtain exact informations for the asymptotics of their coefficients
$[z^N]E_{s}(z)$ and $[z^N]B_{s}(z)$ at leading order for
$N\rightarrow\infty$.  In fact, the asymptotic behaviour of the
coefficients of a generating function is strictly related to the
behaviour of the generating function around its dominant singularity,
i.e.\ its pole of smallest modulus.  For a complete review on
singularity analysis refer to \cite{flajolet2009}, and for a summary
of the basic results that we will need in the following see
Appendix~\ref{app:singanalysis}.

First of all, we need to study how the operator
$\hat{L}_{(\omega_p,\eps)}$ alters the singular behaviour of a
generating function.  Let $\tilde{\bigO}$ be the usual big-O notation,
for expansions of functions of $z$ near $z=1$, but with possible
additional logarithmic factors $\log(1-z)$ that we are not interested
in tracking (that is,
$f(z) = g(z) + \tilde{\bigO}((1-z)^b)$ means that there exists a
finite $c$ such that
$f(z) = g(z) + \bigO((1-z)^b (\ln(1-z))^c)$, which in particular
implies that
$f(z) = g(z) + \bigO((1-z)^{b+\delta})$ for all $\delta>0$)
\begin{proposition}\label{prop:Lq}
    Let $p>0$, $k \in \mathbb{Z}^+$ and let $f(z) = \sum_{N\geq0} f_{N} z^{N}$ be a
    generating function with unit radius of convergence.  
    Suppose that $f$ admits the singular expansion
\begin{equation}
 \begin{split}
     f(z) = f^{\rm reg}(z) + a\, (1-z)^{-\alpha}\left( 1+\tilde\bigO\left( (1-z)^\xi \right)  \right)\quad\text{ for }z\rightarrow1 \, .
 \end{split}
    \end{equation}
with $f^{\rm reg}(z)$ analytic in a neighbourhood of $z=1$ and $\xi > 0$
(in particular, $f(z)$ has a unique dominant singularity, at $z=1$).

Suppose further that 
\begin{equation}
    \begin{split}
\omega_p(N) \sim N^p (1 + \bigO(N^{-\eta})) \quad\text{ for }N\rightarrow\infty \, 
    \end{split}
\end{equation}
for some $\eta >0$.

If 
$\alpha$ and
$\alpha+kp \neq 0,-1,-2,\dots$ then
    \begin{equation}\label{eq:LqGeneral}
 \begin{split}
     \hat{L}_{(\omega_p,\eps)}^k[f](z) &= \tilde{f}^{\rm reg}(z;\omega_p,\eps,f,k) + a \frac{\Gamma\left( \alpha+kp \right)}{\Gamma\left( \alpha \right)} (1-z)^{-(\alpha+kp)} \left( 1+\tilde\bigO\left( \left(1-z \right)^{\min(\eta,\xi)}\right)  \right)\quad\text{ for }z\rightarrow1 \, 
 \end{split}
    \end{equation}
    with a new, \emph{a priori} unknown regular part $\tilde{f}^{\rm reg}(z)$.
\end{proposition}
\pf\
As $\omega_p(N) \sim N^p$ for large $N$, the coefficients of $\hat{L}_{(\omega_p,\eps)}^k[\sum_{N\geq0}f_{N}z^{N}]$ scale as $f_{N}N^{kp}$ for $N\rightarrow\infty$.
The fact that $f_{N}$ gets corrected sub-exponentially (in particular,
algebraically) means that the radius of convergence is not changed by
$\hat{L}_{q}$, so that also $\hat{L}_{q}[f]$ has singular expansion
around $z=1$, and no other singularities of smaller radius.
Now, Corollary~\ref{cor:asy} of Appendix~\ref{app:singanalysis} implies the result at leading order.

An analogous procedure allows to estimate the order of the remainder.
As we are not requiring any condition related to $\xi$ and $\eta$
(namely, that suitable linear combinations are not in a certain range
of integers), there may be additional logarithmic factors, that we
treat consistently by means of the $\tilde\bigO$ notation.  \qed

\bigskip
\noindent

The regular part 
$\tilde{f}^{\rm reg}(z;\omega_p,\eps,f,k)$ is not determined by the
asymptotic behaviour of $f$ and $L(z;\omega_p,\eps)$ alone, and its
calculation requires the full functions $\omega_p$ and $f$, besides
the parameter $\eps$ (see \cite{flajolet2009} for more information on
the singularity analysis of Hadamard products).
In the following, we will be interested in the determination of this
elusive quantity only in the case $k=1$, $\eps=0$ and $f(z)=E_0(z)$,
i.e.\ the regular part at $z=1$ of
\begin{equation}
\label{eq:LhatE0}
    \begin{split}
\hat{L}_{(\omega_p,0)}\left[E_0\right](z) = \sum_{N \geq 0}
c_N \omega_p(N) z^N
=:
\alpha(\omega_p) (1+\bigO(1-z))
+ \frac{\Gamma(p-\frac{1}{2})}{2 \sqrt{\pi}}
(1-z)^{\frac{1}{2}-p}
(1+\bigO((1-z)^{\min\left(\eta,\frac{1}{2}\right)}))
    \end{split}
\end{equation}
where 
$c_k = 2^{-2k-1} C_k =
\frac{\Gamma(k+\frac{1}{2})}{2\sqrt{\pi}\,\Gamma(k+2)}$
are the normalized Catalan numbers (that is, 
$\sum_{k \geq 0} c_k = 1$).

Equation~\eqref{eq:LhatE0}, which corresponds to
(\ref{eq:LhatE0antic}) with a more precise description of the error
terms, implicitly defines the constant $\alpha(\omega_p)$ as the limit
for $z\rightarrow 1$ of the regular part of
$\hat{L}_{(\omega_p,0)}\left[E_0\right](z)$.  For some special choices
of $\omega_p$ it is possible to compute explicitly this regular part,
and we provide some examples in Appendix~\ref{app:costs}.  In any
case, we can formally get rid of the unknown term $\alpha(\omega_p)$
by reabsorbing it into the trivial constant $\epsilon$.

Indeed, for $k=1$, we recall that $\hat{L}_{(\omega_p,\eps)}[f](z) = \hat{L}_{(\omega_p,0)}[f](z)-\eps f(z)$.
Thus, we have that 
\begin{equation}
    \begin{split}
\tilde{f}^{\rm reg}(z;\omega_p,\eps,f,1) = \tilde{f}^{\rm reg}(z;\omega_p,0,f,1) - \eps f^{\rm reg}(z) \, . 
    \end{split}
\end{equation}
As a result we have:
\begin{proposition}
\label{prop:epsilon}
By an appropriate choice of $\eps$, namely
\begin{equation} \label{eq:epsilon}
    \begin{split}
\eps(\omega_p,f) = 
\frac{\tilde{f}^{\rm reg}(1;\omega_p,0,f,1)}{f^{\rm reg}(1)}
\, ,
    \end{split}
\end{equation}
we can guarantee that 
$\tilde{f}^{\rm reg}(z;\omega_p,\eps(\omega_p,f),f,1) = \bigO(1-z)$
for $z \rightarrow 1$.
In the case of $f=E_0$, this reduces to
\be
\epsilon = \alpha(\omega_p)
\ef.
\ee
\end{proposition}
\noindent
We are now ready to study the asymptotic behaviour of the moments in
our models.  By induction, it is easy to prove that the location of
the dominant singularity of $E_{s}(z)$ and of $B_s(z)$ must be at
$z=1$.  Our strategy is to expand Equations~\eqref{eq:recE}~and~\eqref{eq:recB}
to the leading singular order in $1-z$, to obtain
the leading singular order of $E_{s}(z)$ and $B_s(z)$.
\begin{proposition}\label{prop:EBleading}
    Let $\eps = \eps(\omega_p,E_0) = \alpha(\omega_p)$ and $p>0$,
    $p\neq\frac{1}{2}$. 
    Suppose that
    \begin{equation}
\begin{split}
\omega_p(N) \sim N^p (1 + \bigO(N^{-\eta})) \quad\text{ for }N\rightarrow\infty \, 
\end{split}
    \end{equation}
    for some $\eta >0$, and let
    \begin{equation}
\begin{split}
\eta_s =
\begin{cases}
    1 & s=0\\
    \min\left(\eta,\frac{1}{2}\right) & s=1\\
    \min\left(\eta,p,\frac{1}{2}\right) & s\geq2 
\end{cases} \, .
\end{split}
    \end{equation}
Then, for $s\geq 1$
    \begin{equation}\label{eq:ansatzEz}
\begin{split}
E_{s}(z) &= 2 \mu^{\rm (E) }_{s}(p) (1-z)^{-\left(
  p+\frac{1}{2} \right) s+\frac{1}{2}}  
  (1+\tilde\bigO((1-z)^{\eta_s}))
\quad\text{ for } z\rightarrow1\, ,
\end{split}
    \end{equation}
    and
    \begin{equation}\label{eq:ansatzBz}
\begin{split}
B_{s}(z) &= \phantom{2}  \mu^{\rm (B) }_{s}(p)
(1-z)^{-\left( p+\frac{1}{2} \right) s-\frac{1}{2}}
(1+\tilde\bigO((1-z)^{\eta_s}))
  \quad\text{ for } z\rightarrow1\, .
\end{split}
    \end{equation}
    The coefficients $\mu^{\rm (E) }_{s}(p)$ satisfy
    \begin{equation}\label{eq:recMuE}
\begin{split}
\mu^{\rm (E) }_{0}(p) &= -\frac{1}{2} \, , \\
\mu^{\rm (E) }_{1}(p) &= \frac{1}{8\sqrt\pi} \Gamma\left( p-\smfrac{1}{2} \right) \, , \\
\mu^{\rm (E) }_{s}(p) &= \mu^{\rm (E) }_{s-1}(p) \frac{\Gamma\left( s\left( p+\frac{1}{2} \right) -1  \right)}{2\Gamma\left( s\left( p+\frac{1}{2} \right) -p-1 \right)} + \sum_{k=1}^{s-1} \mu^{\rm (E) }_{k}(p)\mu^{\rm (E) }_{s-k}(p) \qquad \text{for} \quad s\geq2
\, .
\end{split}
    \end{equation}
    The coefficients $\mu^{\rm (B)}_{s}(p)$ satisfy
    \begin{equation}\label{eq:recMuB}
\begin{split}
\mu^{\rm (B)}_{0}(p) &= 1\\
\mu^{\rm (B)}_{s}(p) &= 2 \sum_{k=0}^{s-1} \mu^{\rm (B)}_{k}(p)\mu^{\rm (E)}_{s-k}(p) \qquad \text{for} \quad s\geq1 \, .
\end{split}
    \end{equation}
\end{proposition}
\noindent
Note in particular that, if 
$\omega_p(k) = k^p +\bigO(1,k^{p-\frac{1}{2}})$, we can bound
all error-term exponents by $\eta_s = \min(p,\frac{1}{2})$.

A detailed proof can be found in Appendix~\ref{proof:EBleading}.  The
idea is the following.  For $s\in\{0,1\}$,
Equation~\eqref{eq:ansatzEz} and \eqref{eq:ansatzBz}, and the starting
conditions for the recursions in Equation~\eqref{eq:recMuE} and
\eqref{eq:recMuB}, can be verified explicitly.  The particular choice
of $\eps=\alpha(\omega_p)$ is crucial for the ansatz to be correct at
$s=1$, as otherwise the non-null regular part of
$\hat{L}_{(\omega_p,0)}[E_0](z)$ would dominate when $p$ is in the
range $0<p<\frac{1}{2}$.  Then, one proceeds by induction, using
Proposition~\ref{prop:Lq} and the expansions of
Equation~\eqref{eq:recE} and Equation~\eqref{eq:recB} to the leading
singular order in $(1-z)$.

Finally, using again Corollary~\ref{cor:asy} we get that
\begin{equation}
\label{eq:ME}
    \begin{split}
M^{\rm (E)}_s(N;\omega_p,\alpha(\omega_p)) 
  &\sim \frac{4\sqrt{\pi} \, s!}{\Gamma\left( \left(
  p+\frac{1}{2} \right)s -\frac{1}{2} \right)} \mu^{\rm
  (E)}_s(p) N^{s\left( p+\frac{1}{2} \right) }
(1+\tilde\bigO(N^{-\eta_s}))
    \end{split}
\end{equation}
and 
\begin{equation}
\label{eq:MB}
    \begin{split}
M^{\rm (B)}_s(N;\omega_p,\alpha(\omega_p)) 
  &\sim \frac{\sqrt{\pi} \, s!}{\Gamma\left( \left(
  p+\frac{1}{2} \right)s +\frac{1}{2} \right)} \mu^{\rm
  (B)}_s(p) N^{s\left( p+\frac{1}{2} \right) } 
(1+\tilde\bigO(N^{-\eta_s}))
\, .
    \end{split}
\end{equation}
Equations (\ref{eq:ME}) and (\ref{eq:MB}) show that, for both ranges
$p \in (0,\frac{1}{2})$ and $p \in (\frac{1}{2},+\infty)$,
the stochastic leading term in our quantity of interest scales as
$N^{p+\frac{1}{2}}$. Let us define the asymptotic rescaled moments
\begin{equation}
    \wbar{M}_s:=\lim_{N \to \infty} M_s(N;\omega_p,\alpha(\omega_p))\, N^{-s\left( p+\frac{1}{2}
  \right) }
\,,
\end{equation}
that is
\begin{align}
\wbar{M}^{\rm (E)}_s(\omega_p,\alpha(\omega_p)) 
  &= \frac{4\sqrt{\pi} \, s!}
{\Gamma\left( \left( p+\frac{1}{2} \right)s -\frac{1}{2} \right)} 
\mu^{\rm (E)}_s(p)
\\
\wbar{M}^{\rm (B)}_s(\omega_p,\alpha(\omega_p)) 
  &= \frac{\sqrt{\pi} \, s!}
{\Gamma\left( \left( p+\frac{1}{2} \right)s +\frac{1}{2} \right)} 
\mu^{\rm (B)}_s(p)\,,
\end{align}
These are the moments of the candidate distribution our random
variable $x_p$ introduced in Equation~\eqref{eq:introPmaggiore}.
However, we need to prove that these moments, determined by the
recursions in Equation~\eqref{eq:recMuE} and
Equation~\eqref{eq:recMuB}, define \emph{uniquely} two families of
distributions, that we shall call $\rho^{\rm (E/B)}(x;p)$.  
The
uniqueness of these distributions is discussed in the following
section.

As a final remark, let us rephrase the role of the constant
$\epsilon = \alpha(\omega_p)$. Recall that
\begin{equation}
    \begin{split}
A_{(\omega_p)}(w) = \eps N + A_{(\omega_p,\eps)}(w) \, .
    \end{split}
\end{equation}
We proved that the integer moments of $A_{(\omega_p,\eps)}(w)$ scale as $N^{s\left( p+\frac{1}{2} \right) }$, so that, at leading order in $N$ and for $0<p<\frac{1}{2}$, the moments of $A_{(\omega_p)}(w)$ are given by
\begin{equation}
    \begin{split}
\left\langle \left[A_{(\omega_p)}(w)\right]^s \right\rangle \sim \left[ \alpha(\omega_p) \right]^s N^s \, .
    \end{split}
\end{equation}
Thus, for $0<p<\frac{1}{2}$, the distribution of the rescaled variable
$N^{-1} A_{(\omega_p)}$ converges, as $N$ grows, to a Dirac's delta
distribution, with average value $\alpha(\omega_p)$.  In this regime,
the fluctuations around the mean of the rescaled variable are of order
$N^{\frac{1}{2}-p}$, and their distribution, after an appropriate
rescaling, is described by the non-trivial moments
$\overline{M}(\omega_p,\alpha(\omega_p))$.

Notice also that, for $p>\frac{1}{2}$, the candidate distributions are
supported on $[0,+\infty)$ as
they describe a positive random costs, while for $0<p<\frac{1}{2}$ the
distributions are supported on $(-\infty,+\infty)$, which is
compatible with the presence of the shift at leading order.

\subsection{Uniqueness of the distributions $\rho^{\rm (E/B)}(x;p)$}
\label{sec:uniqueness}

The problem of determining whether a moment sequence defines uniquely
a distribution goes under the name of \emph{moment problem}
\cite{lin2017recent}.  In particular, if the distribution is supported
on $[0,+\infty)$, the problem is called 
\emph{Stieltjes moment problem}, while for distributions supported on
$(-\infty,+\infty)$ it is called \emph{Hamburger moment problem}.  In
both cases, a sufficient condition for the uniqueness of the
distribution is given by Carleman's condition.  In the case of the
Hamburger moment problem, the distribution is uniquely determined if
\begin{equation}
    \begin{split}
\sum_{n \geq 1} m_{2n}^{-\frac{1}{2n}} = +\infty \, ,
    \end{split}
\end{equation}
where $\{m_n\}_{n=1}^{+\infty}$ is the moment sequence.
In the case of the Stieltjes moment problem, the distribution is
uniquely determined if
\begin{equation}
    \begin{split}
\sum_{n \geq 1} m_{n}^{-\frac{1}{2n}} = +\infty \, .
    \end{split}
\end{equation}
In both cases, we see that the uniqueness of the distribution
$\rho^{\rm (E/B)}(x;p)$ can be determined by an asymptotic analysis of
$\mu^{\rm (E/B)}_s(p)$ for large order $s$.  If the moment sequences
do not grow too fast, then Carleman's condition will grant the
uniqueness of the distribution.

\begin{proposition}\label{prop:largeS}
    If $p\in \mathbb{R}^+ \smallsetminus \{ \frac{1}{2} \}$, there exist
    $A_p$ and $R_p \in \mathbb{R}^+$ such that
    \begin{equation}
\begin{split}
| \mu^{\rm (E)}_s(p) | &\leq R_p \; A_p^s \; \Gamma\left( ps+1 \right) C_{s-1}\, ,\quad\forall s \geq 1 \\
| \mu^{\rm (B)}_s(p) | &\leq A_p^s \; \Gamma\left( ps+1 \right) C_{s}\, ,\quad\forall s \geq 0 \, .
\end{split}
    \end{equation}
\end{proposition}
\noindent
A proof can be found in Appendix~\ref{proof:largeS}, along with
explicit expressions for $A_p$ and $R_p$. 

By combining Proposition~\ref{prop:largeS} with the normalizations of
Equation~\eqref{eq:ME} and Equation~\eqref{eq:MB}, we find that the
moments $M^{\rm(E/B)}$ grow slower than $\exp\left( \frac{s}{2} \log s
\right)$ for large $s$.  Thus, Carleman's condition immediately
implies that both distributions $\rho^{\rm (E/B)}(x;p)$ are
uniquely determined by their moment sequences for all $p \in
\mathbb{R}^+ \smallsetminus \{\frac{1}{2}\}$.

\section{A combinatorial interpretation of the coefficients 
$\mu^{\rm (E)}_s(p)$}
\label{sec:combinatorial}

\noindent
In this section we show that Equation~\eqref{eq:recMuE} is solved by a
diagrammatic expansion involving rooted planar trees, with a peculiar
form of the weight.

Let us call $\mathcal{T}(s)$ the set of rooted planar trees with $s$
non-root vertices, and, for a tree $T \in \mathcal{T}(s)$ and a vertex
$v \in T$, call $\ell_v$ its out-going degree (i.e., the number of
`children' in the tree), and $k_v$ the number of `descendents' (that
is, w.r.t.\ the notion of hook of a rooted tree introduced in Section
\ref{ssec.hookform}, it is the hook at $v$, minus one).
\begin{proposition}
    \begin{equation}
\begin{split}
\mu^{\rm (E)}_s(p) = \frac{1}{b(s)} \sum_{T \in \mathcal{T}(s)} \prod_{v \in T}
a(\ell_v) b(k_v)  \, 
\end{split}
    \end{equation}
where the range of the product $v \in T$ stands for the $s+1$
vertices of $T$, 
\begin{equation}
\label{eq:Al}
\begin{split}
a(\ell) = \frac{4^{\ell-1}\,\Gamma\left( \ell-\frac{1}{2} \right)}
{\sqrt{\pi}\,\Gamma\left( \ell+1 \right)} \, 
\end{split}
    \end{equation}
    and
    \begin{equation}\label{eq:Bk}
\begin{split}
b(k) = \frac{\Gamma\left( (k+1)(p-\frac{1}{2})+k \right)}
{2\, \Gamma\left( k(p-\frac{1}{2})+k-\frac{1}{2} \right)} \, .
\end{split}
    \end{equation}
\end{proposition}
\pf\,
Given two sequences
$\{a(\ell)\}_{\ell \geq 0}$ and $\{b(k)\}_{k \geq 0}$, we
define the quantities
\begin{equation}
\label{eq:trees}
\begin{split}
\nu(s) = 
\sum_{T \in \mathcal{T}(s)} \prod_{v \in T} 
  b(k_v)a(\ell_v)
\,.
\end{split}
\end{equation}
Let us introduce the generating functions:
\begin{align}
A(z)&=\sum_{\ell} a(\ell) z^{\ell}
\,, \\
\label{eq.defX}
X(z) &= \sum_{s\geq 1} \nu(s-1)
z^s \\
\label{eq.defY}
Y(z) &= \sum_{s \geq 1} 
\frac{\nu(s)}{b(s)}
z^s \, .
\end{align}
Then, the recursive combinatorial definition of rooted planar trees
gives $\nu(0) = b(0)a(0)$ and, for $s \geq 1$,
\begin{equation}
\label{eq:recNu}
\begin{split}
\nu(s) 
&= 
b(s) 
\sum_{\ell \geq 1} 
a(\ell) 
\sum_{\substack{s_1\dots s_\ell \geq 1 \\ \sum_i s_i = s}} 
\prod_{i=1}^\ell \nu(s_i-1) \, ,
\end{split}
\end{equation}
that is, multiplying both sides by $z^s/b(s)$, and summing over 
$s\geq 1$, we get
\be
\label{eq.7856735}
Y(z)
=A(X(z))
-A(0)
\,.
\ee
In the specific case in which $a(\ell)$ are the (shifted) 
Catalan numbers,
\begin{equation}
\label{eq.678765385}
\begin{split}
a(\ell)
= \frac{4^{\ell-1} \Gamma\left( \ell-\frac{1}{2} \right)}{\sqrt{\pi}\,
  \Gamma\left( \ell+1 \right)}
\end{split}
\end{equation}
(i.e.\ $a(0) = -\frac{1}{2}$ and $a(\ell) = C_{\ell-1}$ for
$\ell \geq 1$), we get $A(z) =-\frac{1}{2}\sqrt{1-4z}$, and
Equation~\eqref{eq.7856735} can be rewritten as
\be
\label{eq.7856735b}
Y(z)=X(z)+Y(z)^2
\,.
\ee
(See Table \ref{tab.treeser} for the first few terms of these series.)

By substituting the definitions
(\ref{eq.defX}) and (\ref{eq.defY}), we recognise
Equation~\eqref{eq:recMuE}, under the identifications
\begin{equation}
\begin{split}
b(s) = \frac{\Gamma\left( (s+1)(p-\frac{1}{2}) + s
  \right)}{2\Gamma\left(s(p-\frac{1}{2}) +s-\frac{1}{2} \right)}\,,
\end{split}
\end{equation}
and
$\nu(s) = b(s)\, \mu^{\rm (E)}_s(p)$.
\qed

\begin{table}
\[
\begin{array}{c|c|c|c}
s & \textrm{diagrams} & X(z) & Y(z) \\
\hline
\rule{0pt}{12pt}
0 && b_0 a_0 z & 0
\\
\rule{0pt}{25pt}
1 && b_0 a_0 \, b_1 a_1 \, z^2 & 
b_0 a_0 \, a_1 \, z 
\\
\rule{0pt}{20pt}
&&
&
\\
2 &
\setlength{\unitlength}{10pt}
\begin{picture}(4,1)
\put(0,-.53){\includegraphics{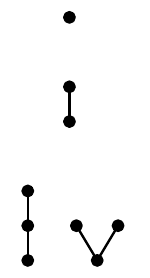}}
\end{picture}
&
\big(b_0 a_0 \, b_1 a_1 \, b_2 a_1 
+ (b_0 a_0)^2 b_2 a_2 \big) z^3
&
\big(b_0 a_0 \, b_1 a_1^2 
+ (b_0 a_0)^2 a_2 \big) z^2
\\
\end{array}
\]
\caption{\label{tab.treeser}\footnotesize The first terms of the series $X(z)$ and
  $Y(z)$ in (\ref{eq.defX}) and (\ref{eq.defY}). Note that Equation
  (\ref{eq.7856735b}) is satisfied up to the order $z^2$, provided
  that $a_1=a_2=1$, which is in agreement with the
  prescription~(\ref{eq.678765385}).}
\end{table}

\section{The case of $p=\frac{1}{2}$}
\label{sec:p05}

\noindent
As already mentioned in the introduction, the moments of 
$\rho^{\rm (E)}(x;p)$, as well as the constant $\alpha(\omega_p)$,
diverge in the limit $p\rightarrow\frac{1}{2}$.  In particular, it is
easy to prove by induction that
\begin{equation}
    \begin{split}
\overline{M}^{\rm (E)}_s(\omega_p,\alpha(\omega_p)) =
\left(\frac{t^*}{ p-\frac{1}{2} }\right)^s 
\left( 1+\bigO\left( p-\smfrac{1}{2} \right)  \right)  \, 
    \end{split}
\end{equation}
where 
\begin{equation}\label{eq:tstar}
    \begin{split}
t^* = 4 \lim_{p\rightarrow\frac{1}{2}} 
\left[
\mu^{\rm (E)}_1(p) \left( p-\smfrac{1}{2} \right) 
\right]
= \frac{1}{2\sqrt{\pi}} \, 
    \end{split}
\end{equation}
is the residue of 
$\overline{M}^{\rm (E)}_1(\omega_p,\alpha(\omega_p))$ at its simple
pole $p=\frac{1}{2}$ (the factor of 4 comes from the prefactor of the
moments).

Moreover, the explicit examples for the computation of
$\alpha(\omega_p)$ presented in Appendix~\ref{app:costs}, plus the
easy universality argument at the end of the same appendix, imply
that, independently from the microscopic details of the function
$\omega_p$,
\begin{equation}
    \begin{split}
\alpha(\omega_p) = -\frac{t^*}{p-\frac{1}{2}} 
\left( 1+\bigO\left( p-\smfrac{1}{2} \right)  \right)
\,.
    \end{split}
\end{equation}
Thus, at least at the level of the first moment, the divergences
cancel out, and the average of $A_{(\omega_p)}$ has a finite limit at
$p=\frac{1}{2}$.

The apparent mechanism is that the two divergences are a spurious
effect of the way we decided to separate the random variable
$A_{(\omega_p)}$ into a deterministic part $\alpha(\omega_p) N$ and a
probabilistic part $x_p N^{p+\frac{1}{2}}$ (see
Equation~\eqref{eq:decomposition1}).  In this case, both
divergences could be regularized by simply shifting the random
variable $x_p N^{p+\frac{1}{2}} \rightarrow \left( x_p - t(p) \right)
N^{p+\frac{1}{2}}$, and reabsorbing the shift into the deterministic
part $\alpha(\omega_p) N \rightarrow \alpha(\omega_p) N + t(p)
N^{p+\frac{1}{2}}$.  The shift function $t(p)$ must satisfy
\begin{equation}
    \begin{split}
t(p) = \frac{t^*}{p-\frac{1}{2}}\left( 1+\bigO\left( p-\smfrac{1}{2}
\right)  \right)  \, .
    \end{split}
\end{equation}
in order to regularize at sight both the deterministic component of
$A_{(\omega_p)}$ and the average value of its probabilistic part. 

If this intuition is right, it remains to be proven that this shift
regularizes all the higher-order moments as well. However, in this
paper we have chosen to use only methods from singularity analysis,
which are badly adapted to non-linear shifts (in $N$), so that all the claims
in the remaining part of this section shall be proven in a second
companion paper in which we perform an asymptotic probabilistic
analysis at fixed $N$ that produces a tree-diagram perturbative
expansion in the spirit of Section \ref{sec:combinatorial} (but with
``coloured'' trees, in order to keep into account the shift terms).

An excerpt of the resulting theory is the following fact:
\begin{definition}
Define the operator 
$\widehat{H}^{\mycirc}_U$, 
acting on functions
$f(\{y_i\}_{i \in U})$, as a multi-dimensional finite-difference
operator:
\be
\widehat{H}^{\mycirc}_U [f(y)] 
:=
\sum_{\{y_j\} \in \{0,1\}^U}
(-1)^{\sum_j y_j}
f(y_1,\ldots,y_{|U|})
\ef,
\ee
and let 
$\widehat{H}_U[f(y)] = 
\lim_{\delta \to 0} 
(-\delta)^{-|U|}
\widehat{H}^{\mycirc}_U[f(\delta y)]$ 
be the corresponding
multi-dimensional L'H\^{o}pital evaluation:
\be
\widehat{H}_U [f] 
:=
\left.
\frac{\dd^{|U|}}{\prod_{i \in U} \dd y_i}
f(y_1,\ldots,y_{|U|})
\right|_{y_1=\cdots=y_{|U|}=0}
\ef.
\ee
\end{definition}

\begin{table}
\[
\begin{array}{|c|c|rcl|}
\hline
\rule{0pt}{12pt}
s & \textrm{diagrams} & 
\multicolumn{3}{|c|}{\textrm{contribution to~}\widetilde{M}_s} \\
\hline
\rule{0pt}{16pt}
1 && 
\left.
\frac{\dd}{\dd y_1} 
\frac{e^{\xi y_1}}{\Gamma(1-y_1)}
a(1) b(1-y_1) 
\right|_{y_1=0}
&= & 0
\\
\rule{0pt}{25pt}
2 && 
\left.
-2
\frac{\dd}{\dd y_1} 
\frac{e^{\xi y_1}}{\Gamma(2-y_1)}
a(1)^2 b(1-y_1) b(2-y_1) 
\right|_{y_1=0}
&= &
\frac{4(\ln 2-1)}{\pi}
\\
\rule{0pt}{10pt}
&
\rule{12pt}{0pt}
&
\rule{12pt}{0pt}
&&
\\
&
\setlength{\unitlength}{10pt}
\begin{picture}(3,1)
\put(0,-.53){\includegraphics{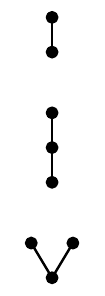}}
\end{picture}
&
\raisebox{-13pt}{\rule{0pt}{3pt}}%
\frac{1}{4 \sqrt{\pi}}
\left.
\frac{\dd^2}{\dd y_1 \dd y_2} 
\frac{e^{\xi (y_1+y_2)}}
{\Gamma(2-y_1-y_2)}
a(2) b(2-y_1-y_2) 
\right|_{y_1=y_2=0}
&= &
\frac{4}{\pi}
-
\frac{\pi}{4}
\\
\hline
\end{array}
\]
\caption{\label{tab.treeser12}\footnotesize The first diagrams
  involved in the evaluation of the moments in Equation
  (\ref{eq.claim12}). Here the functions $a$ and $b$ are as in
  (\ref{eq:Al}) and (\ref{eq:Bk}), evalued at $p=\smfrac{1}{2}$,
  and $\xi=2 \ln 2+\gamma_E$.}
\end{table}

\begin{proposition}
%
For $T$ a rooted tree, calling ${\rm r}$ the root vertex, and $U=U(T)$
the set of leaf vertices. For a vertex $v \in V(T)$ call $U_v$ the set
of leaves which are strictly below $v$.\footnote{In tree order, $u$ is
  strictly below $v$ if $u \neq v$ and the path from ${\rm r}$ to $u$
  goes through $v$.}  Introduce one variable $y_i$ per leaf vertex
$i$, and call
$y_v = \sum_{i \in U_v} y_i$ (in particular, 
$y_{\rm r} = \sum_{i \in U(T)} y_i$).

Under our ``canonical'' choice (\ref{eq.tpCano}) of shift function
$t(p)$, the shifted moments $\widetilde{M}_s$
are given by
\begin{align}
\label{eq.claim12pre}
\widetilde{M}^{\rm (E)}_s(\smfrac{1}{2}+\delta) 
&= 
8 \sqrt{\pi}\, s!
\sum_{T \in \mathcal{T}(s)} 
\widehat{H}^{\mycirc}_{U(T)}
\left[
\frac{
\left(\Gamma(\frac{1}{2})/\Gamma(\frac{1}{2}+\delta)\right)^{y_{\rm r}}
}{
\Gamma\big( s+\delta(s+1-y_{\rm r})\big)
}
\prod_{v \in V(T)}
a(\ell_v) 
b\big(
k_v - \smfrac{\delta}{1+\delta} y_v
\big)
\right]
\,
\end{align}
where $a(\ell)$ and $b(k)$ are as in (\ref{eq:Al}) and
(\ref{eq:Bk}), evalued at $p=\smfrac{1}{2}+\delta$.
In particular, for $\delta=0$,
\be
\begin{split}
\label{eq.claim12}
\widetilde{M}^{\rm (E)}_s(\smfrac{1}{2}) 
&=
8 \sqrt{\pi}\, s!
\sum_{T \in \mathcal{T}(s)} 
(-1)^{|U(T)|}
\widehat{H}_{U(T)}
\left[
\frac{e^{(2 \ln 2+\gamma_E)y_{\rm r}}}
{\Gamma\big( s-y_{\rm r}\big)}
\prod_{v \in V(T)}
a(\ell_v) 
b^{\rm (reg)}(k_v-y_v)
\right]
\\
&=
8 \sqrt{\pi}\, s!
\sum_{T \in \mathcal{T}(s)} 
(-C)^{|U(T)|}
\widehat{H}_{U(T)}
\left[
\frac{e^{(2 \ln 2+\gamma_E)y_{\rm r}}}
{\Gamma\big( s-y_{\rm r}\big)}
\prod_{v \in V(T)\smallsetminus U(T)}
a(\ell_v) 
b(k_v-y_v)
\right]
\end{split}
\ee 
where $b^{\rm (reg)}=b$ for a non-leaf node, and
$C= a(0) b^{\rm (reg)}(0)=\frac{1}{8 \sqrt{\pi}}
=\lim_{\delta \to 0} \delta\, a(0) b(0)$ 
for a leaf node.
\end{proposition}
\noindent
Note that, indeed, 
$\left.\frac{\dd}{\dd p}
\ln(\Gamma(\frac{1}{2})/\Gamma(p))\right|_{p=\frac{1}{2}} 
=2 \ln 2+\gamma_E$,
and that $b(0)$ is the only weight in the expression
(\ref{eq.claim12pre}) that is singular for $\delta \to 0$ (while all
$a(\ell)$ and derivatives $\partial^h b^{\rm (reg)}(k)$ for $h \leq k$
are regular in the limit).

See Table \ref{tab.treeser12} for the first few terms of these series.
As a check of the proposition above, we computed exactly the moments
$\langle (x_p - t(p))^s \rangle$ using Equation~\eqref{eq:recMuEintro}
and the software \emph{Mathematica 12}, with the choice
(\ref{eq.tpCano}) of $t(p)$, up to order $s=5$, finding that they are
all regular at $p=\frac{1}{2}$, and consistent with Equation
(\ref{eq.claim12}) (see Table~\ref{table:regMom}).

\begin{table}
\[
    \begin{array}{|c|c|c|}
\hline 
\rule{0pt}{12pt}\raisebox{-6pt}{\rule{0pt}{12pt}}%
s=2 &
2^3 \ell
-3 \zeta_2
&
0.610375\ldots
\\
\hline 
\rule{0pt}{12pt}\raisebox{-6pt}{\rule{0pt}{12pt}}%
s=3 & 
2^4 \ell(\ell-1)
-8 \zeta_2
+14 \zeta_3
&
0.266217\ldots
\\
\hline
\rule{0pt}{12pt}\raisebox{-6pt}{\rule{0pt}{12pt}}%
s=4 & 
\frac{2^6}{15}
\ell (8 \ell^2+39 \ell+18)
-\frac{2^4}{5} \zeta_2 (61 \ell-8)
+\frac{2^8}{5} \zeta_3
-9 \zeta_2^2
&
1.28827\ldots
\\
\hline
\rule{0pt}{12pt}%
s=5 & 
\frac{2^8}{315}
\ell
(90 \ell^3+1663 \ell^2-1461 \ell-657)
-\frac{2^5}{105}
\zeta_2
(2295 \ell^2+157 \ell+544)
&
1.73555\ldots
\\
\raisebox{-6pt}{\rule{0pt}{12pt}}%
&
+\frac{2^5}{105}
\zeta_3
(5115 \ell-728)
+\frac{4944}{35} \zeta_2^2
-420 \zeta_2 \zeta_3
+744 \zeta_5
&
\\ \hline
    \end{array}
\]
    \caption{\footnotesize Regularized (and rescaled) moments 
$(2 \sqrt{\pi})^s \langle (x_p - t(p))^s \rangle
= (2 \sqrt{\pi})^s \widetilde{M}_s
=2 (2 \sqrt{\pi})^{s+1} 
\frac{\Gamma(s+1)}{\Gamma(s-\frac{1}{2})}
\tilde{\mu}_s$ for 
$t(p)$ as in (\ref{eq.tpCano}), at $p=\frac{1}{2}$. Here $\ell$ is a
shortcut for $\ln 2$, and $\zeta_n$ denotes the Riemann's zeta
function $\zeta(n)$. Note that, under the formal replacement $\zeta_n
\to \ell^n$, the expressions above are rational polynomials in $\ell$,
of degree $s$.}
    \label{table:regMom}
\end{table}

\section{The limit of $p \rightarrow 0$, and logarithmic cost
  functions}
\label{sec:p0log}

\noindent
It is interesting to study the limiting behaviour of $A_{(\omega_p)}$
when $p$ goes to zero along a family of functions 
$\{\omega_p\}_{p \in \mathbb{R}^+}$.  As discussed in the
Introduction, as the limit itself is trivially 
$A_{(\omega_p)}(w) = N$, we shall instead focus on the behaviour of
the rescaled shifted quantity
\begin{equation}\label{eq:p0}
    \begin{split}
B_{(\omega_p)}(w) = \frac{A_{(\omega_p)}(w) - N}{p} \, 
    \end{split}
\end{equation}
as $p \rightarrow 0^+$.  Notice that this is equivalent to the study
of the original observable $A_{(\omega_{0^+})}$, where by
$\omega_{0^+}$ we denote a positive regular function such that
$\omega_{0^+}(k) \sim \log(k)$ for large $k$.  The strategy of Section
\ref{sec:proofNo05} remains valid, but, as our treatment was assuming that the
function $\omega$ has a power-law asymptotics, in this section we have
to repeat our analysis.
For simplicity of notation, in the remainder of this
section we drop the $+$ superscript of $\omega_{0^+}$.

In the case of Dyck excursions, we will prove that, in the limit of
large $N$,
\begin{equation}
    \begin{split}
A_{(\omega_0)} \das \alpha(\omega_0) N + \sqrt{\gamma_E N \log N} z  +
\bigO\big( \sqrt{N} \big)  \,  
    \end{split}
\end{equation}
where $z$ is a standard Gaussian random variable, $\gamma_E$ is the
Euler's gamma constant and $\alpha(\omega_0)$ is a non-universal
constant that depends on the details of $\omega_0(k)$ away from the
large $k$ regime (see Equation~(\ref{eq:LhatE0})).
To this end, we will compute the integer moments of the observable
$A_{(\omega_0,\epsilon)}$, where $\epsilon$ will be tuned to directly
access the fluctuations around the mean, and we will find that they
are equal to the moments of a Gaussian random variable.
The case of Dyck bridges can be
  treated analogously.

We follow very closely the line of thought of Section~\ref{sec:proofNo05}, taking for given the definitions and the results presented therein.
In particular, Proposition~\ref{prop:recursionGF} holds in the logarithmic case as well.
Thus, we will start our treatment of the problem by generalizing Proposition~\ref{prop:Lq} to the case of logarithmic functions $\omega_0$.
In the following, we define 
\begin{equation}
    \begin{split}
L(z) = \log\left( \frac{1}{1-z} \right)  \, .
    \end{split}
\end{equation}

\begin{proposition}\label{prop:LqLOG}
    Let $k \in \mathbb{Z}^+$ and $f(z) = \sum_{N\geq0} f_{N} z^{N}$ be a
    generating function with unit radius of convergence.  
    Suppose that $f$ admits the singular expansion
\begin{equation}
 \begin{split}
     f(z) = f^{\rm reg}(z) + a\, (1-z)^{-\alpha} \left[ \log\left( \frac{1}{1-z} \right)\right]^\beta  \left( 1+\bigO\left( \left[L(z)\right]^{-1} \right)  \right)\quad\text{ for }z\rightarrow1 \, .
 \end{split}
    \end{equation}
with $f^{\rm reg}(z)$ analytic in a neighbourhood of $z=1$ and $\beta >0$
(in particular, $f(z)$ has a unique dominant singularity, at $z=1$).

Suppose further that 
\begin{equation}
    \begin{split}
\omega_0(N) \sim (\log N) (1 + \bigO(N^{-\eta})) \quad\text{ for }N\rightarrow\infty \, 
    \end{split}
\end{equation}
for some $\eta >0$.

If 
$\alpha \neq 0,-1,-2,\dots$ then
    \begin{equation}
 \begin{split}
     \hat{L}_{(\omega_0,\eps)}^k[f](z) &= \tilde{f}^{\rm reg}(z;\omega_0,\eps,f,k) + a (1-z)^{-\alpha} \left[ L(z) \right]^{\beta + k}  \left( 1+\bigO\left( \left[ L(z)\right]^{-1} \right)  \right)\quad\text{ for }z\rightarrow1 \, 
 \end{split}
    \end{equation}
    with a new, \emph{a priori} unknown regular part $\tilde{f}^{\rm reg}(z)$.
    If, instead, $\alpha = 0$, then
    \begin{equation}
 \begin{split}
     \hat{L}_{(\omega_0,\eps)}^k[f](z) &= \tilde{f}^{\rm reg}(z;\omega_0,\eps,f,k) + a  \left[ L(z) \right]^{\beta + k+1}  \left( 1+\bigO\left( \left[ L(z)\right]^{-1} \right)  \right)\quad\text{ for }z\rightarrow1 \, .
 \end{split}
    \end{equation}
\end{proposition}
\pf\
See the proof of Proposition~\ref{prop:Lq}.
Here one must pay attention to logarithmic factors; see \cite[Sec. VI.2]{flajolet2009}.
Notice that the final error term here is independent on the error terms on $f(z)$ and on $\omega_0(k)$ as long as they are algebraic.
\qed

Now that we have specified the leading behaviour of
$\hat{L}_{(\omega_0, \epsilon)}$, we can again study the singular
behaviour of $E_s(z)$ inductively.
\begin{proposition}\label{prop:EleadingLOG}
    Let $\eps = \eps(\omega_0,E_0) = \alpha(\omega_0)$.
    Suppose that
    \begin{equation}
\begin{split}
\omega_0(N) \sim (\log N) (1 + \bigO(N^{-\eta})) \quad\text{ for }N\rightarrow\infty \, 
\end{split}
    \end{equation}
    for some $\eta >0$.
    Then,
    \begin{equation}\label{eq:ansatzEzLog}
\begin{split}
E_{s}(z) &=
\begin{cases}
-\frac{1}{2} L(z) \left( 1+\bigO\left( [L(z)]^{-1} \right)  \right)  & s=1 \\
    \tau_{s} (1-z)^{\frac{1-s}{2}} \left[ L(z)  \right]^{\frac{s}{2}}
+\bigO\left( (1-z)^{\frac{1-s}{2}}[L(z)]^{\frac{s-2}{2}} \right) & s\geq2
\end{cases}
\end{split}
    \end{equation}
    for $z \rightarrow 1$.
The coefficients $\tau_{s}$ satisfy 
$\tau_2=\frac{\gamma_E}{4}$, and
    \begin{equation}\label{eq:recTauE}
\begin{split}
\tau_s =
\begin{cases}
\displaystyle{
C_{\ell-1} 
2^{1-\ell} 
\tau_2^\ell
}
   & s = 2 \ell, \quad \ell\geq2\\
    0      & s = 2 \ell +1, \quad \ell \geq 1 \\
\end{cases}
\, .
\end{split}
    \end{equation}
Implicitly, $\tau_1=0$, because $E_1(z) \sim L(z)$ instead of 
$L(z)^{\frac{1}{2}}$.
\end{proposition}

\pf\, We prove the result by induction.
The expression for $E_1(z)$
can be easily obtained by using Equation~\eqref{eq:recE} and
Proposition~\ref{prop:LqLOG}.  Here it is crucial to set $\epsilon =
\alpha(\omega_0)$ to discard a contribution to $E_1(z)$ of order
$(1-z)^\frac{1}{2}$ that would otherwise dominate the entire induction
process.

For $s=2$, to obtain the leading singular order of $E_2(z)$, one must
expand the term $\hat{L}_{(\omega_0,\epsilon)}[E_1](z)$ to its
next-to-leading order.  This verifies that the scaling given in
Equation~\eqref{eq:ansatzEzLog} and the value of $\tau_2$ are correct.

Finally, for $s\geq3$, it is easy to verify, using
Equation~\eqref{eq:recE}, that the leading singular order of $E_s(z)$
has the correct scaling, and that

\begin{equation}
\label{eq.398756874}
    \begin{split}
\tau_s = \frac{1}{2} \sum_{k=2}^{s-2} \tau_k \tau_{s-k} \, .
    \end{split}
\end{equation}
This recursion is consistent with the ansatz that $\tau_k=0$ for all
odd $k$, which is indeed implied by the fact that $\tau_1=0$.
Moreover, by introducing the generating function
\begin{equation}
    \begin{split}
T(z) = \sum_{s=1}^\infty \tau_{2s}\;z^s \, 
    \end{split}
\end{equation}
it is easy to recognise in (\ref{eq.398756874}) the classical
quadratic relation for Catalan numbers, and in turns
see that
\begin{equation}
    \begin{split}
T(z) = 1-\sqrt{1-2\tau_2 z}  \, ,
    \end{split}
\end{equation}
whose expansion gives the explicit expression for $\tau_{2s}$ given in Equation~\eqref{eq:recTauE}.
\qed

Thus, by using Corollary~\ref{cor:asy}, we find that the  even integer moments $M_{2\ell}^{\rm (E)}(N;\omega_0, \alpha(\omega_0))$ satisfy
\begin{equation}
    \begin{split}
M_{2 \ell}^{\rm (E)}(N; \omega_0, \alpha(\omega_0)) \sim (2\ell-1)!! ( \gamma_E N \log N)^{\ell}
    \end{split}
\end{equation}
where $s!!$ denotes the double factorial, 
and that the odd integer moments are null at order $(N \log N)^\frac{s}{2}$.
Thus, the rescaled moments $\overline{M}_s^{\rm (E)}$ satisfy
\begin{equation}
    \begin{split}
\overline{M}_{2 \ell}^{\rm (E)}(\omega_0, \alpha(\omega_0)) 
= \lim_{N\rightarrow \infty} 
\frac{M_{2 \ell}^{\rm (E)}(N; \omega_0, \alpha(\omega_0))}{(N \log N)^\ell}
= \gamma_E^\ell (2\ell-1)!!
 \, 
    \end{split}
\end{equation}
which are precisely the integer moments of a Gaussian random variable
with variance $\gamma_E$.

We finally prove Equation~(\ref{eq:p0}) by observing that:
\begin{itemize}
    \item the $\epsilon = \alpha(\omega_0)$ shift accounts for the deterministic term of Equation~\eqref{eq:p0};
    \item Proposition~\ref{prop:EleadingLOG}, along with Carleman's condition, accounts for the stochastic term of Equation~\eqref{eq:p0};
    \item the dominant error term is not induced by the error on the higher moments, but by the leading behaviour of the average value, that scales as $\sqrt{N}$.
\end{itemize}
The computations of this section can be easily generalized to Dyck
bridges, and to cost functions with asymptotic behaviour $(\log k)^p$
with $p>0$, but we do not detail this here.

\appendix

\section{Some facts about the area-Airy distribution}
\label{app:airyfacts}

\noindent
We recall the main known facts about the area-Airy distribution, following
the review \cite{Majumdar2005}. The area-Airy distribution 
$f(x) = f_{\rm Ai}(x)$ has support on $\mathbb{R}^+$, and Laplace
transform \cite{darling1983,louchard1984}
\begin{equation}
    \begin{split}
\tilde{f}_{\rm Ai}(\lam) = \int_{0}^{\infty} f_{\rm Ai}(x)
e^{-\lam x} dx 
= \lam \sqrt{2\pi}\sum_{k=1}^{\infty}e^{-a_{k} \lam^{2/3} 2^{-1/3}}
 \, 
    \end{split}
\end{equation}
where $a_{k}$ is the value of the $k$-th zero of the standard Airy function $\Ai(x)$.
A closed expression for the density was found in \cite{takacs1991}
\begin{equation}
    \begin{split}
f_{\rm Ai}(x) = \frac{2\sqrt{6}}{x^{10/3}} \sum_{k=1}^{\infty} e^{-b_{k}/x^{2}}b_{k}^{2/3}U\left( -\frac{5}{6},\frac{4}{3},\frac{b_{k}}{x^{2}} \right) 
\, 
    \end{split}
\end{equation}
where $b_{k}=2a_{k}^{3}/27$ and $U(a,b,z)$ is the confluent hypergeometric function. 
Moreover, a recursion for the moments of $f_{\rm Ai}(x)$ is known
\cite{takacs1991}: define $K_s$ as 
\begin{equation}
    \begin{split}
M^{\rm Ai}_{s} = \int_{0}^{\infty} f_{\rm Ai}(x) x^{s} dx 
=: \sqrt{\pi}\,2^{(4-s)/2}
\frac{\Gamma\left( s+1 \right)}{\Gamma\left( \frac{3s-1}{2} \right)} K_{s}\, .
    \end{split}
\label{eq.2876547}
\end{equation}
Then
\begin{equation}
\label{eq:recTakacs}
    \begin{split}
K_{s} &= \frac{3s-4}{4} K_{s-1} + \sum_{j=1}^{s-1} K_{j}K_{s-j} \quad \forall s \geq 1 \\
K_{0} &= -\frac{1}{2} \,.
\, 
    \end{split}
\end{equation}
Finally, $f_{\rm Ai}(x)$ has asymptotic behaviours
\cite{takacs1991,takacs1995}
\begin{equation}
    \begin{array}{rlp{2cm}l}
f_{\rm Ai}(x) &\sim x^{-5} e^{-2a_{1}^{3}/27x^{2}} & & \text{ as } x\rightarrow0\\
f_{\rm Ai}(x) &\sim  e^{-6x^{2}}& & \text{ as } x\rightarrow\infty\,.
    \end{array}
\end{equation}

\section{Some facts about singularity analysis}
\label{app:singanalysis}

\noindent
The main result that we will need is the following theorem (here
stated informally, see \cite{flajolet2009} for a precise statement):
\begin{theorem}\label{thm:asy}
    Let $f(z) = \sum_{N\geq0}f_{N}z^{N}$ and $g(z) =
    \sum_{N\geq0}g_{N}z^{N}$ be two generating functions with radius
    of convergence $r$. Then
\begin{equation}
 \begin{split}
     f(z) \sim g(z) \quad\text{ for }z\rightarrow r \quad\iff\quad
     f_{N} \sim g_{N} 
\quad\text{ for }N\rightarrow\infty \, .
 \end{split}
\end{equation} 
\end{theorem}
\noindent
In particular, we will need two special cases:
\begin{corollary}
\label{cor:asy}
    Let $f(z) = \sum_{N\geq0}f_{N}z^{N}$ be a generating function with
    unit radius of convergence. If $f$ admits the singular expansion
    \begin{equation}
\label{eq.29857435647}
 \begin{split}
     f(z) = f^{\rm reg}(z) + a\, 
(1-z)^{-\alpha} \left( \log\frac{1}{1-z} \right)^{\beta} 
\left( 1+\mathcal{O}\left( 1-z \right)  \right)
\quad\text{ for }z\rightarrow1 \, 
 \end{split}
    \end{equation}
    with $\alpha \neq 0,-1,-2,\dots$ and $f^{\rm reg}(z)$ analytic in a neighbourhood of $z=1$, then
    \begin{equation} \label{eq:thmasy}
 \begin{split}
     f_{N} 
=
\frac{a}{\Gamma\left( \alpha
       \right)}N^{\alpha-1}\left( \log N \right)^{\beta}  
(1+\bigO(N^{-1}))
\sim 
\frac{a}{\Gamma\left( \alpha
       \right)}N^{\alpha-1}\left( \log N \right)^{\beta}  \quad\text{
       for } N\rightarrow\infty\, ,
 \end{split}
    \end{equation} 
    and the viceversa is also true.
    If $\alpha=0$, $\beta=1$ the same statement holds with
    \begin{equation}
 \begin{split}
     f_{N} \sim \frac{1}{N}
      \, .
 \end{split}
    \end{equation}
\end{corollary}
\pf\
Follows from the expansions given in \cite{flajolet2009} for products
of polynomial and logarithmic singularities.
\qed\\

Similar statements hold if the error term has a different
form. Namely, if in (\ref{eq.29857435647}) we replace
$\mathcal{O}\left( 1-z \right)$ by
$\mathcal{O}\Big( 
(1-z)^{-a} \big( \log\frac{1}{1-z} \big)^{b}
\Big)$, 
in (\ref{eq:thmasy}) we get an error term of the form
$1+\bigO(N^{-a}(\ln N)^b)$ if $\alpha+a\neq0,-1,-2,\dots$
Cases in in which $\alpha+a=0,-1,-2,\dots$ must be treated separately; see \cite{flajolet2009}.

\section{Examples of cost functions}
\label{app:costs}

\noindent
In this Appendix, we provide some examples of families of cost
functions $\omega_p$ for which we can compute analytically the
quantity $\alpha(\omega_p)$ defined in (\ref{eq:LhatE0}), that is the
regular part at $z=1$ of
\begin{equation}
    \begin{split}
\hat{L}_{(\omega_p,0)}\left[E_0\right](z) = \sum_{N \geq 0} c_N \omega_p(N) z^N \,.
    \end{split}
\end{equation}
As a result, in these cases it is possible to compute exactly the
shift $\eps(\omega_p,E_0)$, that just coincides with
$\alpha(\omega_p)$, and have a complete control over the asymptotic
behaviour of the random variable $A_{(\omega_p)}$ given in
(\ref{eq:decomposition1}). Recall that $c_N = 2^{-2N-1} C_N$ are the
normalized Catalan numbers, and that the regular part at $z=1$ of
$E_0(z)$ equals~$1$.

As a first example, we set 
$\omega^{(\frac{1}{2})}_p(k) = 
\frac{\Gamma\left( k+p+\frac{1}{2} \right)}{\Gamma\left( k+\frac{1}{2} \right)}$.
In this case
\begin{equation}
\begin{split}
    \hat{L}_{(\omega^{(\frac{1}{2})}_p,0)}\left[E_0\right](z) &=
    \frac{\Gamma\left( p-\frac{1}{2} \right)}{2\sqrt{\pi} \, z}
    \left( \left( 1-z \right)^{\frac{1}{2}-p} -1  \right) \, 
\end{split}
\end{equation}
whose regular part at $z=1$ equals
$-\frac{\Gamma\left( p-\frac{1}{2} \right)}{2\sqrt{\pi} \, z} $, giving
\begin{equation}
\begin{split}
\alpha(\omega^{(\frac{1}{2})}_p) = \eps(\omega^{(\frac{1}{2})}_p,E_0) 
= -\frac{\Gamma\left( p-\frac{1}{2} \right)}{2\sqrt{\pi}} 
= \frac{\Gamma\left( p-\frac{1}{2} \right)}{\Gamma\left( -\frac{1}{2} \right)}
\, .
\end{split}
\end{equation}
As a second example, we set $\omega^{(1)}_p(k) = \frac{\Gamma\left( k+p+1 \right)}{\Gamma\left( k+1 \right)}$.
In this case
\begin{equation}
\begin{split}
\hat{L}_{(\omega^{(1)}_p,0)}\left[E_0 \right](z) = \frac{\Gamma\left( p+1 \right)}{2} \pFq{2}{1}{\frac{1}{2},1+p}{2}{z} 
\end{split}
\end{equation}
where $_2F_1$ is the hypergeometric function,
whose regular part at $z=1$ can be extracted by using the inversion
formula \cite[Equation~15.3.6, pg.~559]{abramowitz1972},
giving
\begin{equation}
\begin{split}
\alpha(\omega^{(1)}_p) = -\frac{
\Gamma\left( p \right) \Gamma\left( \frac{1}{2}-p \right)
}{\Gamma\left( \frac{1}{2} \right)\Gamma\left(-p\right)}  \, .
\end{split}
\end{equation}
More generally, for
$\omega^{(a)}_p(k) = \frac{\Gamma\left( k+p+a \right)}{\Gamma\left(
  k+a \right)}$, 
defined for $a+p$ not a negative integer, and positive
for all integer $k$ when $a>0$, we get
\begin{equation}
\begin{split}
\hat{L}_{(\omega^{(a)}_p,0)}\left[E_0 \right](z) = 
\frac{\Gamma\left( a+p-1 \right)}
{z \Gamma\left( a-1 \right)}
\left(
1 - \pFq{2}{1}{-\frac{1}{2},a+p-1}{a-1}{z} 
\right)
\,.
\end{split}
\end{equation}
Again, the regular part at $z=1$ of the hypergeometric function can be
extracted by using the inversion formula, giving
\begin{equation}
\begin{split}
\alpha(\omega^{(a)}_p) = 
\Gamma\left( a+p-1 \right)
\left(
\frac{1}{\Gamma\left( a-1 \right)} - 
\frac{\Gamma\left( \frac{1}{2}-p \right)}
{\Gamma\left( a-\frac{1}{2} \right)\Gamma\left( -p \right)}
\right)
\end{split}
\end{equation}
Note that, within this family, only the case $a=\frac{1}{2}$ gives an
expression for $\alpha(\omega^{(a)}_p)$ which is finite for all 
$p \in \mathbb{R}^+ \smallsetminus \{\frac{1}{2}\}$.

Another simple family is
$\omega^{(a,-\frac{3}{2})}_p(k) = 
\frac{\Gamma\left( k+p+a-\frac{3}{2} \right)
\Gamma\left( k+2 \right)}
{\Gamma\left(k+a \right)
\Gamma\left( k+\frac{1}{2} \right)
}$, 
defined for $a+p-\frac{3}{2}$ not a negative integer, and generalising 
$\omega^{(\frac{1}{2})}$ (as 
$\omega^{(2,-\frac{3}{2})}_p(k) = \omega^{(\frac{1}{2})}_p(k)$),
that gives
\begin{equation}
\begin{split}
\hat{L}_{(\omega^{(a,-\frac{3}{2})}_p,0)}\left[E_0 \right](z) = 
\frac{\Gamma\left( a+p-\frac{3}{2} \right)}
{2 \sqrt{\pi}\; \Gamma\left( a \right)}
\pFq{2}{1}{1,a+p-\frac{3}{2}}{a}{z} 
\,,
\end{split}
\end{equation}
and in turns
\begin{equation}
\begin{split}
\alpha(\omega^{(a,-\frac{3}{2})}_p) = 
\frac{
\Gamma\left( a+p-\frac{3}{2} \right)
}
{
\Gamma\left( \frac{1}{2} \right)
\Gamma\left( a-1 \right)
(1-2p)
}
\,.
\end{split}
\end{equation}
This quantity is finite for all 
$p \in \mathbb{R}^+ \smallsetminus \{\frac{1}{2}\}$ whenever
$a>\frac{3}{2}$.

Finally, we observe that, in all the cases analysed here in detail, we
have that, in the limit $p \to \frac{1}{2}$,
\be
\alpha(\omega_p)
=
\frac{1}{p-\frac{1}{2}}
\left(
-\frac{1}{2\sqrt{\pi}} 
+
o(1)
\right)
\ef.
\ee
This is not a coincidence. Indeed, the operator $\hat{L}_{(\omega,0)}$
is linear: 
$\hat{L}_{(\omega'+\omega'',0)} f(z)
=\hat{L}_{(\omega',0)} f(z)
+\hat{L}_{(\omega'',0)} f(z)$,
so that, if $\omega'$ is any of the families above, and $\omega''$ is
a family of functions of the form 
$\omega''_p(k)=k^p (1+b(k,p) k^{-\eta})$, with $b(k,p)$ uniformly
bounded and $\eta>0$, then 
$(\omega'_p-\omega''_p)(k)=k^{p-\eta} b'(k,p)$, with $b'(k,p)$
uniformly bounded, and
$\lim_{z \to 1}
|\hat{L}_{(\omega'_{\frac{1}{2}}-\omega''_{\frac{1}{2}},0)} E(z)|
< +\infty$. As the singular parts of 
$\hat{L}_{(\omega'_{\frac{1}{2}},0)} E(z)$ and of
$\hat{L}_{(\omega''_{\frac{1}{2}},0)} E(z)$ are the same
(because they are determined only by $p$), and are diverging, it must
be the case that also the regular (in $z$) parts diverge (in $p$) with
the same coefficient.

\section{Proof of Proposition~\ref{prop:recursionGF}}
\label{proof:recursionGF}

\noindent
By using the fact that a Dyck excursion $w$ can always be
decomposed as $w = (+,w_{2},-,w_{1})$, with
$w_{1}$ and $w_{2}$ being Dyck excursions of length $m$ and
$N-m-1$ (for some $0 \leq m \leq N-1$) and $+/-$ an up/down step, one finds that

\begin{equation}
\begin{split}
 &\frac{1}{s!} \sum_{w \in \mathcal{C}_{N}}\left[
    A_{(\omega_p,\eps)}(w) \right]^{s} 
\\ 
&= \frac{1}{s!}
  \sum_{m=0}^{N-1} 
\sum_{w_{2} \in \mathcal{C}_{m}}
\sum_{w_{1} \in \mathcal{C}_{N-m-1}} 
\left[ A_{(\omega_p,\eps)}(w_{1}) +
    A_{(\omega_p,\eps)}(w_{2}) + \left( \omega_p(m) - \eps\right)
    \right]^{s} 
\\ 
&= \frac{1}{s!} 
  \sum_{m=0}^{N-1} 
\sum_{w_{2} \in \mathcal{C}_{m}}
\sum_{w_{1} \in \mathcal{C}_{N-m-1}} 
  \sum_{\substack{s_{1},s_{2},s_{3}\geq 0 \\ s_{1}+s_{2}+s_{3}=s}}
  \left[ A_{(\omega_p,\eps)}(w_{1})\right]^{s_{1}}\left[
    A_{(\omega_p,\eps)}(w_{2})\right]^{s_{2}}\left[\omega_p(m)-\eps\right]^{s_{3}}
  \frac{s!}{s_{1}!s_{2}!s_{3}!}
\\ 
&= \sum_{m=0}^{N-1}
  \sum_{\substack{s_{1},s_{2},s_{3}\geq 0 \\ s_{1}+s_{2}+s_{3}=s}}
  \left[ \frac{1}{s_{1}!}\sum_{w \in C_{N-m-1}} \left[
      A_{(\omega_p,\eps)}(w) \right]^{s_{1}} \right]
  \frac{\left[\omega_p(m)-\eps\right]^{s_{3}} }{s_{3}!} 
  \left[\frac{1}{s_{2}!} \sum_{w \in C_{m}} \left[
      A_{(\omega_p,\eps)}(w) \right]^{s_{2}} \right]
\end{split}
\end{equation}
implying that 
\begin{equation}
\begin{split}
 2 E_{s}(z) = \delta_{s,0} + z  \sum_{\substack{s_{1},s_{2},s_{3}\geq0
     \\ s_{1}+s_{2}+s_{3}=s}} E_{s_{1}}(z) \frac{1}{s_{3}!}
 \hat{L}_{(\omega_p,\eps)}^{s_{3}}[E_{s_{2}}](z)
   \, ,\qquad s\geq0 \, .
\end{split}
\label{eq.287657442}
\end{equation}
The case $s=0$ gives rise to an equation involving only $E_0(z)$, and
no summations, which is nothing but the generating function of
normalized Catalan numbers $c_k := 2^{-2k-1} C_k$.
So we easily get
\begin{equation}\label{eq:E0}
\begin{split}
 E_{0}(z) = \frac{1-\sqrt{1-z}}{z} \, ,
\end{split}
\end{equation}
For the general case $s \geq 1$, notice that the term $E_{s}(z)$
appears in both sides of the $s$-th Equation~(\ref{eq.287657442}), and only
linearly.  If we isolate $E_{s}(z)$, we obtain after some simplifications
\begin{equation}
\begin{split}
E_{s}(z) = \frac{z}{2\sqrt{1-z}}
\sum_{\substack{s_{1},s_{2},s_{3}\geq0
    \\ s_{1}+s_{2}+s_{3}=s\\s_{1},s_{2}<s}} E_{s_{1}}(z)
\frac{1}{s_{3}!} \hat{L}_{(\omega_p,\eps)}^{s_3}[E_{s_{2}}](z)
\, . 
\end{split}
\end{equation}
The proof for Dyck bridges is in the same spirit of the one for Dyck
paths.  Now, a Dyck bridge decomposes uniquely as 
$\pm w=(+,w_{2},-,w_{1})$, where $w_{2}$ is a Dyck excursion and
$w_{1}$ is a Dyck bridge, hence the recursion involving $E_{s_{2}}$.
Note that, as the first step of a Dyck bridge can be either a $+$ or a
$-$, a factor of two must be taken into account when dealing with this
decomposition.

\section{Proof of Proposition~\ref{prop:EBleading}}
\label{proof:EBleading}

\noindent
We give the proof for $E_s(z)$, the one for $B_s(z)$ being completely
analogous, and we proceed by induction.

The case $s=0$ is trivially verified from the explicit form of
$E_0(z)$.  

The case $s=1$ can be computed explicitly by using
Equation~\eqref{eq:recE}:
\begin{equation}
\begin{split}
E_1(z) 
&= \frac{z}{2\sqrt{1-z}} E_0 \hat{L}_{(\omega_p,\eps)}[E_0](z) =
\frac{1-\sqrt{1-z}}{2\sqrt{1-z}} \hat{L}_{(\omega_p,\eps)}[E_0](z) \,
.
\end{split}
\end{equation}
Using Proposition~\ref{prop:Lq} with $k=1$, $f(z) = E_0(z)$ and thus $f^{\rm
  reg}(z) = \frac{1}{z}$, $\alpha=-\frac{1}{2}$ and $\xi=1$ we obtain
\begin{equation}
\begin{split}
\hat{L}_{(\omega_p,\eps)}[E_0](z) &= \tilde{E}_0^{\rm
  reg}(z;\omega_p,\eps,E_0,1) - \frac{\Gamma\left( p-\frac{1}{2}
  \right)}{-2\sqrt{\pi}}(1-z)^{\frac{1}{2}-p}\left( 1+\tilde\bigO(
  (1-z)^{\min(\eta,1)})
\right) \\
  &= \frac{\Gamma\left( p-\frac{1}{2}
  \right)}{2\sqrt{\pi}}(1-z)^{\frac{1}{2}-p}
  \left( 1+\tilde\bigO\left( \left( 1-z \right)^{\min(\eta,1,\frac{1}{2}+p)} \right)  \right) 
\end{split}
\end{equation}
where the second passage is due to the choice $\eps =
\eps(\omega_p,E_0)$ as in Equation~\eqref{eq:epsilon} to eliminate the
regular part at $z=1$.  Thus,
\begin{equation}
\begin{split}
E_1(z) &= \frac{1-\sqrt{1-z}}{2\sqrt{1-z}} 
 \frac{\Gamma\left( p-\frac{1}{2} \right)}
{2\sqrt{\pi}}(1-z)^{\frac{1}{2}-p} 
\left( 1+\tilde\bigO\left( \left( 1-z \right)^{\min(\eta,1,\frac{1}{2}+p)} \right)  \right) 
\\ &=
\frac{\Gamma\left( p-\frac{1}{2} \right)}{4 \sqrt{\pi}}(1-z)^{-p} 
\left( 1+
\tilde\bigO\left((1-z)^{\min(\eta,\frac{1}{2})}\right) 
\right)
\end{split}
\end{equation}
(where the exponents $1$ and $\frac{1}{2}+p$ in the error term 
are always subleading w.r.t.\ the exponent $\frac{1}{2}$ coming from
the algebraic prefactor).

Notice that if we had left $\eps$ free, the leading signularity of
$E_1(z)$ would have had exponent $-\frac{1}{2}$ for $p<\frac{1}{2}$.
Thus, the tuning of $\eps$ is crucial to allow for a unified treatment
for all $p\neq\frac{1}{2}$.

For $s\geq2$, we suppose that
Equation~\eqref{eq:ansatzEz} is correct for $E_m(z)$, $0\leq m\leq
s-1$, with an error term of the form $1+\tilde\bigO((1-z)^{\eta_m})$, and
we compute the singular expansion around $z=1$ of
Equation~\eqref{eq:recE}.  First of all, Proposition~\ref{prop:Lq}
tells us that, for $s_3 \geq 1$,
\begin{equation}
\begin{split}
\hat{L}^{s_{3}}_{(\omega_p,\eps)}\left[ E_{s_{2}}(z) \right] =
\tilde{E}_{s_{2}}^{\rm reg}(z) + \frac{2 \mu^{\rm (E)}_{s_{2}}(p)
  \Gamma\left( \left(p+\frac{1}{2}\right)s_{2}-\frac{1}{2}+ps_{3}
  \right)}{\Gamma\left( \left(p+\frac{1}{2}\right)s_{2}-\frac{1}{2}
  \right)}(1-z)^{-\left(p+\frac{1}{2}\right)s_{2}+\frac{1}{2}-ps_{3}}
  (1+\tilde\bigO((1-z)^{\min(\eta,\eta_{s_2})}))
\,
\end{split}
\end{equation}
while for $s_3=0$ the specialisation of the RHS to this value holds,
with the simpler error term $1+\bigO((1-z)^{\eta_{s_2}})$.

Notice that all the non-integrity conditions for the singular
exponents  of Proposition~\ref{prop:Lq} are satisfied under our ansatz because $s_{2}<s$ and
$s_{3}\geq0$, and that $\tilde{E}_{s_2}^{\rm reg}(z)$ has various
dependences, that we drop for simplicity.
Then, Equation~\eqref{eq:recE} reduces to
\begin{equation}
\label{eq:recEsingularleading}
\begin{split}
 &2 \mu^{\rm (E)}_{s}(p) (1-z)^{1-\left(p+\frac{1}{2}\right)s}\left( 1+\tilde\bigO\left( \left( 1-z \right)^{\eta_s}  \right)  \right) = \\
 &\quad=\sum_{\substack{s_{1},s_{2},s_{3}\geq0\\s_{1}+s_{2}+s_{3}=s\\s_{1},s_{2}<s}} \frac{1}{s_{3}!} \left[
 E_{s_{1}}^{\rm reg}(z)\tilde{E}_{s_{2}}^{\rm reg}(z) \right.   +\\
 &\quad+ \tilde{E}_{s_{2}}^{\rm reg}(z) \mu^{\rm (E)}_{s_{1}}(p) (1-z)^{\frac{1}{2}-\left(p+\frac{1}{2}\right)s_{1}} \left( 1+\bigO\left( \left( 1-z \right)^{\eta_{s_1}} \right)  \right)+ \\
 &\quad +
 E_{s_{1}}^{\rm reg}(z) \frac{\mu^{\rm (E)}_{s_{2}}(p)\Gamma\left( \left(p+\frac{1}{2}\right)s_{2}-\frac{1}{2}+ps_{3} \right)}{\Gamma\left( \left(p+\frac{1}{2}\right)s_{2}-\frac{1}{2} \right)} (1-z)^{ \frac{1}{2}-\left(p+\frac{1}{2}\right)s_{2}-ps_{3}} 
 \left( 1+\tilde\bigO\left( \left( 1-z \right)^{\min(\eta_{s_2},\eta)} \right)  \right)+ \\
 &\left.\quad+ 
2 \frac{\mu^{\rm (E)}_{s_{1}}(p)\mu^{\rm (E)}_{s_{2}}(p)\Gamma\left( \left(p+\frac{1}{2}\right)s_{2}-\frac{1}{2}+ps_{3} \right)}{\Gamma\left( \left(p+\frac{1}{2}\right)s_{2}-\frac{1}{2} \right)}
 (1-z)^{1-\left(p+\frac{1}{2}\right)(s_{1}+s_{2})-ps_{3}} \left( 1+\tilde\bigO\left( \left( 1-z \right)^{\min(\eta_{s_1},\eta_{s_2},\eta)} \right)  \right)\right]
  \, 
\end{split}
\end{equation}
where $E_{s_{1}}^{\rm reg}(z)$ is the regular part of $E_{s_{1}}(z)$,
and the equality is at leading order in powers of $(1-z)$.

In the RHS of Equation~\eqref{eq:recEsingularleading}, only the third
and fourth term contribute to the leading order, and the former only
for $(s_1,s_2,s_3)=(0,s-1,1)$, while the latter only for
$(s_1,s_2,s_3)=(k,s-k,0)$, with 
$1 \leq k \leq s-1$; this immediately
gives the recursion in Equation~\eqref{eq:recMuE} for the 
$\mu^{\rm (E)}_{s}(p)$ coefficients.

It is also easy to verify that our claim on the form of the error-term
exponents $\eta_s$ holds inductively, by explicitly analysing the
subleading contributions of the various terms in
Equation~\eqref{eq:recEsingularleading}.


\section{Proof of Proposition~\ref{prop:largeS}}
\label{proof:largeS}

First of all, we study the case of Dyck excursions. We drop the
superscript (E) and the dependence on $p$ for simplicity.
Equation~\eqref{eq:recMuE} implies that
\begin{equation}
\begin{split}
|\mu_s| \leq \frac{\Gamma\left( s\left( p+\frac{1}{2} \right) -1
  \right)}{2\Gamma\left( s\left( p+\frac{1}{2} \right) -p-1 \right)}
|\mu_{s-1}| + \sum_{k=1}^{s-1} |\mu_k||\mu_{s-k}| \, .
\end{split}
\end{equation}
We want to prove by induction that
\begin{equation}
\label{eq:ansatzLargeS}
\begin{split}
| \mu_s | \leq R\, A^s\, \Gamma\left( ps+1 \right) C_{s-1}\qquad\forall s \geq 1 
\end{split}
\end{equation}
for some values of $R$ and $A$ (possibly depending on $p$).
For $s=1$, we obtain that $R$ and $A$ must satisfy the condition
\begin{equation}
\begin{split}
RA &\geq f(p)
\end{split}
\end{equation}
where
\begin{equation}
\begin{split}
f(p)&= 
\frac{\left| \Gamma\left( p-\frac{1}{2} \right)\right|}
{8\sqrt{\pi}\,\Gamma\left(p+1  \right) }\, .
\end{split}
\end{equation}
For the inductive step, we assume that the ansatz in
Equation~\eqref{eq:ansatzLargeS} is satisfied for all $\mu_k$ with 
$1 \leq k \leq s-1$ for some fixed values of $R$ and $A$.  This
implies, using Equation~\eqref{eq:recMuE}, that
\begin{equation}\label{eq:cond3}
\begin{split}
| \mu_s | &\leq 
R\, A^s\, \Gamma\left( ps+1 \right) C_{s-1}
\left[ \frac{\Gamma\left( s\left( p+\frac{1}{2} \right) -1  \right)}{2\Gamma\left( s\left( p+\frac{1}{2} \right) -p-1 \right)}\frac{1}{A}\frac{\Gamma\left( p(s-1)+1 \right)}{\Gamma\left( ps+1 \right)} \frac{C_{s-2}}{C_{s-1}} + \right. \\
  &\hspace{3.75cm}\left. 
+ R \sum_{k=1}^{s-1} \frac{C_{k-1}C_{s-k-1}}{C_{s-1}}\frac{\Gamma\left( pk+1 \right)\Gamma\left( p(s-k)+1 \right)}{\Gamma\left( ps+1 \right)}
\right] \, .
\end{split}
\end{equation}

The first term in the parenthesis of Equation~\eqref{eq:cond3} must be treated separately for different values of $s$:
\begin{itemize}
\item for $s=2$, it equals exactly $\frac{1}{4A}$;
\item for $s$ even and $s\geq 4$, it equals
\begin{equation}\label{eq:footnote}
\begin{split}
\frac{1}{2A} \frac{s}{4s-6} \prod_{i=1}^{\lfloor s/2 \rfloor-2}
\frac{sp+\lfloor s/2 \rfloor-i-1}{(s-1)p+\lfloor s/2 \rfloor-i-1} \leq
\frac{1}{2A}\frac{s}{4s-6}\left( \frac{sp+1}{(s-1)p+1}
\right)^\frac{s-4}{2} \, 
\end{split}
\end{equation}
where the last term is
bounded from above by its limit for $s,p\rightarrow\infty$,%
\footnote{The supremum over $p$ of the function at the right-hand side of
  Equation~\eqref{eq:footnote} is realised for 
$p\rightarrow \infty$, 
and equals
\begin{equation}
\begin{split}
\frac{s}{4s-6}\left( \frac{s}{s-1} \right)^{\frac{s-4}{2}}  \, .
\end{split}
\end{equation}
This function changes monotonicity at the zeroes with odd multiplicity
of the function
\begin{equation*}
\begin{split}
-3+(2s-3)(s-1) \left( s \log\left( \frac{s}{s-1} \right)  - 1 \right)
\, .
\end{split}
\end{equation*}
From the fact that $\log\left( \frac{s}{s-1} \right) \leq \frac{1}{s} + \frac{1}{2s^2}$ for $s\geq2$, 
we get that there are no zeroes on the right of the right-most zero of the equation
$6s = (2s-3)(s-1)$,
which is slightly smaller than 6. Thus, the only candidate
maxima are $s=4$ and the limit $s \to +\infty$, and the
latter is larger than the former.}
   that is $\frac{\sqrt{e}}{8A}$;
\item for $s=3$, it is monotone decreasing in $p$, and can be
  similarly bounded by its limit for $p\rightarrow 0$, i.e.\ $\frac{1}{4A}$;
\item for $s$ odd and $s\geq 5$, it can be bounded by
\begin{equation}
\begin{split}
   \frac{1}{2A}\frac{s}{4s-6}\frac{\Gamma\left( s\left( p+\frac{1}{2} \right) -\frac{1}{2}  \right)}{\Gamma\left( s\left( p+\frac{1}{2} \right) -p-\frac{1}{2} \right)}\frac{\Gamma\left( p(s-1)+1 \right)}{\Gamma\left( ps+1 \right)} =
\frac{1}{2A}\frac{s}{4s-6} \prod_{i=1}^{\lfloor s/2 \rfloor-2}
\frac{sp+\lfloor s/2 \rfloor-i-1}{(s-1)p+\lfloor s/2 \rfloor-i-1}
   \, 
\end{split}
\end{equation}
that, in turn, can be bounded by the same procedure used in the even
case (which, incidentally, produces the same value for the bound).
\end{itemize}
The resulting bound is the maximum over the bounds of the various
terms. As $\sqrt{e}<2$, we have that the first term of
Equation~\eqref{eq:cond3} can be bounded from above by $\frac{1}{4A}$,
for all $s\geq2$ and all $p>0$.

The second term in the parentheses of Equation~\eqref{eq:cond3} can be
simplified by using the fact that the Gamma function is
logarithmically convex, giving in particular that
\begin{equation}
\begin{split}
\log \Gamma\left( t x +1 \right) \leq t \log \Gamma\left( x+1 \right)\, ,\quad\forall t \in [0,1] \, .
\end{split}
\end{equation}
Thus,
\begin{equation}
\begin{split}
\frac{\Gamma\left( pk+1 \right)\Gamma\left( p(s-k)+1 \right)}{\Gamma\left( ps+1 \right)} \leq 
\frac{\Gamma\left( ps+1 \right)^\frac{k}{s}\Gamma\left( ps+1 \right)^\frac{s-k}{s}}{\Gamma\left( ps+1 \right)} = 1
 \, .
\end{split}
\end{equation}
At this point, the sum over $k$ gives $1$, thanks to the well-known
recursion for Catalan numbers, i.e.
\begin{equation}
\begin{split}
\sum_{i=0}^{n-1} C_i C_{n-i-1} = C_n \, .
\end{split}
\end{equation}
As a result, we have obtained two condition that must be satisfied by
$R$ and $A$ to confirm that the ansatz in
Equation~\eqref{eq:ansatzLargeS} is indeed true:
\begin{equation}\label{eq:cond1}
\begin{split}
RA \geq f(p)
\end{split}
\end{equation}
and
\begin{equation}\label{eq:cond2}
\begin{split}
\frac{1}{4A}+R \leq 1 \, .
\end{split}
\end{equation}
For the case of Dyck bridges, Equation~\eqref{eq:recMuB} implies
\begin{equation}
\label{eq.4987663575}
\begin{split}
|\mu^{\rm (B)}_s| \leq  2\sum_{k=0}^{s-1} |\mu^{\rm (B)}_k||\mu^{\rm (E)}_{s-k}| \, .
\end{split}
\end{equation}
We now have the ansatz
\begin{equation}
\begin{split}
|\mu^{\rm (B)}_s| \leq A^s \Gamma\left( sp+1 \right)C_{s} \, , \quad \forall s \geq 0
 \, 
\end{split}
\end{equation}
where $A$ is the same as for Dyck excursions. Substituting the
equation above, and (\ref{eq:ansatzLargeS}), into
(\ref{eq.4987663575}), easily shows that the ansatz
is verified if $R$ satisfies
\begin{equation}
\label{eq:cond3a}
\begin{split}
R &\geq \frac{1}{2} \, .
\end{split}
\end{equation}
One choice of functions $R_p$ and $A_p$ that satisfies all three
conditions, Equations~\eqref{eq:cond1},~\eqref{eq:cond2} and~\eqref{eq:cond3a} is given by
\begin{equation}
\begin{split}
(A_p,R_p)
 &= 
\begin{cases}
\left(
\frac{1}{2},
\frac{1}{2}
\right)
& \text{if}\quad 0 < f(p) \leq \frac{1}{4}\\
\left(
\frac{1+4f(p)}{4},
\frac{4f(p)}{1+4f(p)}
\right)
& \text{if}\quad f(p) > \frac{1}{4}
\end{cases} \, ,
\end{split}
\end{equation}
Notice that Equation~\eqref{eq:cond1} cannot be satisfied for
$p=\frac{1}{2}$, as $f(p)$ (and thus also $A_p$) diverge in this limit.

\bibliographystyle{unsrt}
\bibliography{Moments}

\end{document}